\documentclass{article}
\usepackage{amsmath}
\usepackage{amsfonts}
\usepackage{amssymb}
\usepackage{color}
\usepackage{graphicx}
\usepackage{wasysym}
\usepackage{epsfig}
\usepackage{hyperref}
\usepackage{multicol}
\usepackage[utf8]{inputenc}
\usepackage{authblk}
\usepackage[titletoc]{appendix}
\usepackage{tikz}
\usetikzlibrary{positioning}

\usepackage{float}
\usepackage{caption}
\newfloat{Digraph}{htbp}{loa}
\newfloat{SP}{htbp}{sp}[Digraph]

\newtheorem{theorem}{Theorem}

\newtheorem{example}[theorem]{Example}

\def\qed{\hfill $\Box$\medskip}

\title{On Sign Pattern Matrices that Allow or Require Algebraic Positivity}

\author[1]{Jean Leonardo Abagat \thanks{jlabagat@gmail.com}}
\author[1,2]{Diane Christine Pelejo \thanks{dcpelejo@math.upd.edu.ph}}

\affil[1]{Institute of Mathematics,  College of Science, University of the Philippines Diliman}
\affil[2]{ Natural Sciences Research Institute,  College of Science, University of the Philippines Diliman}

\date{}
\begin{document}
\openup .5\jot
\maketitle

\begin{abstract}

A square matrix $M$ with real entries is said to be \textit{algebraically positive} (AP) if there exists a real polynomial $p$ such that all entries of the matrix $p(M)>0$. A square sign pattern matrix $S$ is said to \textit{allow} algebraic positivity if there is an algebraically positive matrix $M$ whose sign pattern class is $S$. On the other hand, $S$ is said to \textit{require} algebraic positivity if any matrix $M$, having sign pattern class $S$, is algebraically positive. Motivated by open problems raised in the work of Kirkland, Qiao and Zhan (2016) on AP matrices, we list down all nonequivalent irreducible $3\times 3$ sign pattern matrices and classify each of them into three groups (i) those that require AP, (ii) those that allow but not require AP, or (iii) those that do not allow AP. We also give a necessary condition for an irreducible $n\times n$ sign pattern to allow algebraic positivity.
\end{abstract}


\section{Introduction}
We denote the set of $m\times n$ real matrices by $M_{m\times n}$ and the set of all $n\times n$ matrices by $M_n$ and the $(i,j)$ entry of a matrix $A$ by $(A)_{ij}$. 

A matrix $A\in M_{m\times n}$ is said to be \textit{positive}, written $A>0$, if all entries of $A$ are positive. Several applications of positive matrices appear in various fields of study such as Markov chains in probability theory, population models, iterative methods in numerical analysis, economic models, epidomiology, low-dimensional topology, physics, and many more (for example, see\cite{MacCluer,Pillai}). Some natural generalizations of the concept of matrix positivity have also been studied extensively such as matrix nonnegativity, matrix primitivity and eventual positivity. In \cite{Kirk1}, Kirkland, Qiao and Zhan presented another generalization called algebraic positivity. A matrix $A\in M_n$ is said to be \textit{algebraically positive} if $p(A)>0$ for some real polynomial $p$. Using the Cayley-Hamilton Theorem, it can be deduced that a matrix $A\in M_n$, where $n\geq 2$, is algebraically positive if and only if there are real numbers $k_1,\ldots, k_{n-1}$ such that the off-diagonal entries of the following matrix are all positive
\[k_1A+\cdots +k_{n-1}A^{n-1}.\] 
The authors of \cite{Kirk1} used the Perron-Frobenius Theorem to prove the following characterization of algebraically positive matrices. 

\begin{theorem}[Kirkland, Qiao and Zhan, 2016]\label{thm1}
A real matrix is algebraically positive if and only if it has a simple real eigenvalue and corresponding left and right positive eigenvectors.
\end{theorem} 

A \textit{sign pattern} is a matrix whose entries are from the set $\{0,+,-\}$. We define the \textit{pattern class} of a sign pattern $S$, denoted by $Q(S)$, as the set
\[Q(S)=\{A\in M_{m,n}\ | \ \mbox{sgn}((A)_{ij})=(S)_{ij}\}\]
of real matrices whose sign pattern is $S$. In this paper, we are interested in determining when a given square sign pattern matrix \textit{allows algebraic positivity}, i.e. there exists an algebraically positive matrix in its sign pattern class; or \textit{requires algebraic positivity}, i.e. all matrices in its sign pattern class are algebraically positive.  

The following results from \cite{Kirk1} will be useful in our study:
\begin{theorem}[Kirkland, Qiao and Zhan, 2016]\label{thm2} \phantom{ \ } \ 
\begin{enumerate}
\item Every algebraically positive matrix is irreducible.
\item If $A$ is an irreducible real matrix all of whose off-diagonal entries are non-negative (or nonpositive), then $A$ is algebraically positive.
\item If a sign pattern allows algebraic positivity, then every row and column contains a $\ +$, or every row and column contains a $\ - $. 
\end{enumerate}
\end{theorem}

Recall that a matrix $A\in M_n$ is reducible if there exists a permutation matrix $P$ such that $PAP^T$ is of the form 
\[PAP^T=\begin{bmatrix}
A_1 & A_2\\
0 & A_3
\end{bmatrix}\]
for some $A_1\in M_s$ such that $1\leq s<n$. The reducibility of a matrix is easily determined by looking at its associated digraph. The \textit{digraph} of $A$, denoted by $\Gamma(A)$, is the directed graph whose vertex set is $\{1,\ldots,n\}$ and whose edge set is $E_A=\{(i,j) \ | \ (A)_{ij}\neq 0\}$. It is known (see \cite{handbook}) that a matrix is irreducible if and only if its digraph is strongly connected, i.e. for every pair of distinct vertices $i$ and $j$, there is a directed path in $\Gamma(A)$ from $i$ to $j$. Note that the digraph of $A^T$ is obtained from the digraph of $A$ by reversing the direction of all the edges. Meanwhile, if $P$ is a permutation matrix, then the digraph of $PAP^T$ is obtained from the digraph of $A$ by permuting the labels of the vertices of $\Gamma(A)$. Thus, if $\Gamma(A)$ is irreducible, then so are $\Gamma(A^T)$ and $\Gamma(PAP^T)$. For the purpose of this study, we say that two digraphs are \textit{equivalent digraphs} if one can be obtained from the other by (i) reversing the direction of all the edges, (ii) permuting the labels of the vertices, or (iii) doing both. In Appendix A, we list down all 26 nonequivalent irreducible digraphs in with vertex set $\{1,2,3\}$.

\begin{theorem}\label{thm3}
If $A$ is algebraically positive, then the following matrices are also algebraically positive 
\begin{multicols}{2}
\begin{enumerate}
\item $A^T$;
\item $-A$;
\item $PAP^T$ for any permutation matrix $P$;
\item $\beta A + \alpha I$ for any $\alpha\in\mathbb{R}$ and any $\beta\in \mathbb{R}\setminus\{0\}$.
\end{enumerate}
\end{multicols}
\end{theorem}
\textit{Proof:} Suppose $p$ is a real polynomial such that $p(A)>0$ and suppose $P$ is a permutation matrix. Then $p(A^T)=p(A)^T>0$ and $p(PAP^T)=Pp(A)P^T>0$. Define $q(x)=p(\frac{1}{\beta}x-\frac{\alpha}{\beta})$, then $q(\beta A + \alpha I)=p(A)>0$. In particular, if $\beta=-1$ and $\alpha=0$, we have $q(-A)>0$. \qed

For the purpose of our study, we say that two sign pattern matrices are \textit{equivalent sign patterns} if one can be obtained from the other by (i) transposition; (ii) permutation similarity; (iii) negation; (iv) any combination of the first three transformations. 

\section{Main Results}
The main result of this study gives a necessary condition for sign patterns that allow algebraic positivity. 
\begin{theorem} \label{thm4}
	Suppose A is an $n \times n$ sign pattern matrix. We write $A = A_+ + A_-$ where $A_+$ contains only the nonnegative signs of $A$, while $A_-$ contains only non-positive signs of $A$. Define $B_A = A_+ - A_-^T$. If $A$ is irreducible and $B_A$ is reducible, then $A$ does not allow algebraic positivity.
\end{theorem}
\textit{Proof:} Let $A\in M_n(\{+,-,0\})$ and
\[(A_+)_{ij}=\left\{\begin{array}{ll}
+ & \mbox{ if }(A)_{ij}=+\\
0 & \mbox{ otherwise}\\
\end{array}\right. \quad 
(A_-)_{ij}=\left\{\begin{array}{ll}
- & \mbox{ if }(A)_{ij}=-\\
0 & \mbox{ otherwise}\\
\end{array}\right. \]
Define $B_A=A = A_+ - A_-^T$. Clearly $B_A\in M_n(\{0,+\})$. Suppose $B_A$ is reducible. Then, we can assume without loss of generality that 
\begin{equation}\label{redB}
	B_A = \begin{bmatrix}
	B_1 & B_2 \\
	0   & B_3
	\end{bmatrix}
\end{equation}
for some $B_1\in M_s(\{0,+\})$, $s \times (n-s)$ sign pattern matrix $B_2$, and $(n-s) \times (n-s)$ sign pattern matrix $B_3$. Otherwise, we can replace $A$ by $PAP^T$ wherer $P$ is a permutation matrix such that $PB_AP^T$ is in the said form. From ($\ref{redB}$), it will follow that $A$ is of the form \[A = \begin{bmatrix}
	C - D & E \\
	-F     & G - H
	\end{bmatrix}\]
for some $C,D\in M_s(\{0,+\})$, $G,H\in M_{n-s}(\{0,+\})$ and $E\in M_{s,n-s}(\{0,+\})$ and $F\in M_{n-s,s}(\{0,+\})$. Note also that the entries of $E$ or $F$ cannot all be zero since $A$ is irreducible. 
	
	Suppose there exists an $n \times n$ matrix $X \in Q(A)$ such that $X$ is algebraically positive. Then we can say that \[X = \begin{bmatrix}
	X_C - X_D & X_E \\
	-X_F       & X_G - X_H
	\end{bmatrix}\in Q(A)\]
	where $X_J\in Q(J)$ for $J=C,D,E,F,G,H$. By Theorem 1, $X$ has a simple real eigenvalue, say $\lambda \in \mathbb{R}$, with corresponding positive left and right eigenvectors \[\begin{bmatrix}
	u^T & v^T \\
	\end{bmatrix}, \begin{bmatrix}
	y\\
	z
	\end{bmatrix} > 0, \mbox{ where } u, y \in M_{s \times 1} \mbox{ and } v, z \in M_{(n-s) \times 1}.\]	Note that 
	\[\begin{bmatrix}
	u^T & v^T \\
	\end{bmatrix} (X - \lambda I) = \begin{bmatrix}
	u^T(X_C - X_D - \lambda I) - v^TX_F & \ast 
	\end{bmatrix}
	= 
	\begin{bmatrix}
	0 & 0 & \cdots & 0
	\end{bmatrix}\]
	and \[(X - \lambda I) \begin{bmatrix}
	y\\
	z
	\end{bmatrix}
	=
	\begin{bmatrix}
	(X_C - X_D - \lambda I)y + X_Ez \\
	\\
	\ast
	\end{bmatrix}
	=
	\begin{bmatrix}
	0 \\ 0 \\ \vdots \\ 0
	\end{bmatrix}\]
	Therefore, $u^T(X_C - X_D - \lambda I) = v^TX_F$ and $ (X_C - X_D - \lambda I)y = -X_Ez$.	Consider the matrix $u^T(X_C - X_D - \lambda I)y = v^TX_Fy = -u^TX_Ez$. Since the entries of $u^T,v^T, y,z$ are all positive and $X_E,X_F$ are nonzero nonnegative matrices, then $v^TX_Fy$ is a positive real number, while $-u^TX_Ez$ is a negative real number. This gives a contradiction. Hence, any $X\in Q(A)$ is not algebraically positive. \qed\\
Here is an example on how to utilize this result.
\begin{example}\label{ex4.2}
	\normalfont Consider the sign pattern $A$ given below. We write $A$ as follows, separating the positive ($A_+$) and the negative $(A_-)$ signs.
	\[A= \begin{bmatrix}
	0 & + & 0 \\
	+ & 0 & - \\
	+ & 0 & +
	\end{bmatrix} = \begin{bmatrix}
	0 & + & 0 \\
	+ & 0 & 0 \\
	+ & 0 & +
	\end{bmatrix} + \begin{bmatrix}
	0 & 0 & 0 \\
	0 & 0 & - \\
	0 & 0 & 0
	\end{bmatrix}\]
	We then write its $B_A$ as follows, negating and transposing $A_-$.
	$$B= \begin{bmatrix}
	0 & + & 0 \\
	+ & 0 & 0 \\
	+ & + & +
	\end{bmatrix} = \begin{bmatrix}
	0 & + & 0 \\
	+ & 0 & 0 \\
	+ & 0 & +
	\end{bmatrix} + \begin{bmatrix}
	0 & 0 & 0 \\
	0 & 0 & 0 \\
	0 & + & 0
	\end{bmatrix}$$

We can treat the negative signs of $A$ as a reversal of the direction on the corresponding edges in the signed digraph of $A$ to obtain the signed digraph of $B_A$ as seen below.
\bigskip\\
\begin{minipage}{0.49\textwidth}
	\centering
	\begin{tikzpicture}[->,>=stealth,shorten >=1pt,auto,node distance=2cm,
	thick,main node/.style={circle,draw,font=\bfseries}]
	
	\node[main node] (1) {1};
	\node[main node] (2) [below left of=1] {2};
	\node[main node] (3) [below right of=1] {3};
	
	\path
	(1) edge [transform canvas={xshift=-2pt,yshift=2pt},shorten <= -1pt] node[above]{$+$} (2)
	
	(2) edge [transform canvas={xshift=2pt,yshift=-2pt},shorten <= -1pt] node[below]{$+$} (1)
	edge [transform canvas={xshift=0pt,yshift=-3pt},shorten <= -1pt] node[below]{$-$} (3)
	
	(3) edge [loop below] node{$+$} (3)
	edge [transform canvas={xshift=2pt,yshift=2pt},shorten <= -1pt] node[above]{$+$} (1); 
	\end{tikzpicture}
	\\Signed Digraph of $A$
\end{minipage}
\begin{minipage}{0.49\textwidth}
	\centering
	\begin{tikzpicture}[->,>=stealth,shorten >=1pt,auto,node distance=2cm,
	thick,main node/.style={circle,draw,font=\bfseries}]
	
	\node[main node] (1) {1};
	\node[main node] (2) [below left of=1] {2};
	\node[main node] (3) [below right of=1] {3};
	
	\path
	(1) edge [transform canvas={xshift=-2pt,yshift=2pt},shorten <= -1pt] node[above]{$+$} (2)
	
	(2) edge [transform canvas={xshift=2pt,yshift=-2pt},shorten <= -1pt] node[below]{$+$} (1)
	
	(3) edge [loop below] node{$+$} (3)
	edge [transform canvas={xshift=2pt,yshift=2pt},shorten <= -1pt] node[above]{$+$} (1)
	edge [transform canvas={xshift=0pt,yshift=3pt},shorten <= -1pt] node[above]{$+$} (2);  
	\end{tikzpicture}
	\\Signed Digraph of $B_A$
\end{minipage}\bigskip\\
Note that $B_A$ is reducible since its corresponding digraph s not strongly connected \cite{handbook}.
\end{example}
We also note that the converse of Theorem 3 is not true, by providing a counterexample below.
\begin{example}
	\normalfont Consider the sign pattern $C$ below. By the Theorem 2.3, $C$ does not allow algebraic positivity. Now, we write $C$ as follows:
	$$C = \begin{bmatrix}
	0 & - & 0 \\
	- & 0 & + \\
	+ & 0 & +	
	\end{bmatrix} = 
	\begin{bmatrix}
	0 & 0 & 0 \\
	0 & 0 & + \\
	+ & 0 & +	
	\end{bmatrix} + 
	\begin{bmatrix}
	0 & - & 0 \\
	- & 0 & 0 \\
	0 & 0 & 0	
	\end{bmatrix}$$
	However, one can easily check that the matrix $B_C$ below is irreducible.
	\[B_C = \begin{bmatrix}
	0 & + & 0 \\
	+ & 0 & + \\
	+ & 0 & +	
	\end{bmatrix} = 
	\begin{bmatrix}
	0 & 0 & 0 \\
	0 & 0 & + \\
	+ & 0 & +	
	\end{bmatrix} + 
	\begin{bmatrix}
	0 & + & 0 \\
	+ & 0 & 0 \\
	0 & 0 & 0	
	\end{bmatrix}.\]
\end{example}
Recall that a matrix $X\in M_n$ is algebraically positive if and only if $X+\alpha I_n$ is algebraically positive for any $\alpha\in \mathbb{R}$. We now introduce the notion of a \textit{scalar-shift subclass of a sign pattern matrix}. We say that a sign pattern $B$ is a (scalar-shift) subclass of another sign pattern $A$, denoted by $B\trianglelefteq A$, if for any $X\in Q(B)$, there exists $\alpha\in \mathbb{R}$ such that the sign pattern of $X+\alpha I_n$ is equivalent to $A$.
\begin{example}
Let $B=\begin{bmatrix}
0 & +\\
+ & 0
\end{bmatrix}$ and $A=\begin{bmatrix}
- & + \\
+ & -
\end{bmatrix}$. Then $B\trianglelefteq A$ but $A\not \trianglelefteq B$. 
\end{example}
The following theorem follows directly from Theorem 3.4. 
\begin{theorem} Suppose $A \trianglelefteq  B$.
\begin{enumerate}
\item If $B$ requires algebraic positivity, then so does $A$.
\item If $A$ allows algebraic positivity, then so does $B$. 
\end{enumerate}
\end{theorem}
Finally, in Appendix A, we list down all nonequivalent irreducible $3\times 3$ sign pattern matrices and determine whether they require, allow or do not allow algebraic positivity.


\section{Acknowledgments}
This work is funded by the Natural Sciences Research Institute of the College of Science, University of the Philippines Diliman under the research grant MAT-18-1-03. 

\nocite*
%

\newpage
\appendix
\section{Nonequivalent Irreducible $3\times 3$ Sign Pattern Matrices}
We will say that a given sign pattern matrix is (a) RAP if it requires algebraic positivity; (b) AAP if it does not require but allows algebraic positivity; (c) DNA if it does not allow algebraic positivity. 

In this section, we wish to classify each possible $3 \times 3$ sign pattern matrix into one of the three types above. Since reducible matrices belong to the DNA category, we will only consider irreducible matrices. Moreover, we will only list down nonequivalent sign pattern matrices. For efficiency, we list the nonequivalent digraphs and group them  according to their edge count. Then for each digraph, we enumerate all corresponding nonequivalent sign pattern matrices. 

\renewcommand{\theenumi}{\Alph{enumi}}
\renewcommand{\theenumii}{\arabic{enumii}}
\renewcommand{\theenumiii}{\alph{enumiii}}
\begin{enumerate}
	\item \textbf{Graphs with 3 directed edges}. Up to equivalence, there is only one strongly connected digraph with 3 edges and there are 2 nonequivalent sign pattern matrices having this digraph. We classify these sign patterns using Theorems \ref{thm2}.2 and \ref{thm2}.3 to, as indicated below. 
	\medskip\\
	\begin{minipage}{0.4\textwidth}
		\centering
		\includegraphics[scale=0.6]{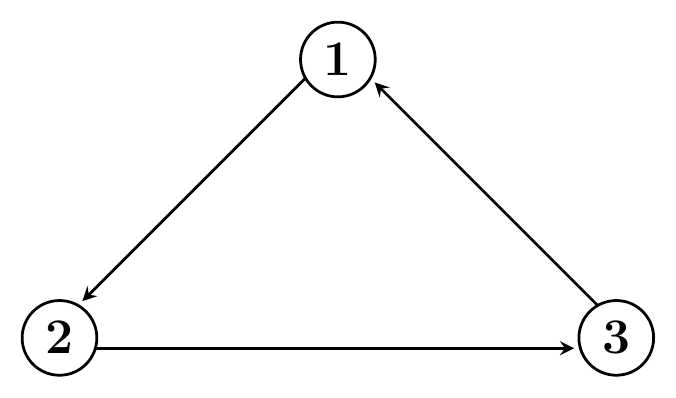}
		\captionof{Digraph}{}
	\end{minipage}
	\begin{minipage}{0.28\textwidth}
		\centering 
		$\begin{bmatrix} 
		0 & + & 0 \\ 
		0 & 0 & + \\
		+ & 0 & 0 
		\end{bmatrix}$
		\captionof{SP}{\textbf{RAP}}
	\end{minipage}
	\begin{minipage}{0.28\textwidth}
		\centering
		$\begin{bmatrix} 
		0 & - & 0 \\
		0 & 0 & + \\ 
		+ & 0 & 0 
		\end{bmatrix}$
		\captionof{SP}{\textbf{DNA}}
	\end{minipage}
	\item \textbf{Graphs with 4 directed edges}. There are three nonequivalent strongly connected digraphs with 4 directed edges. We list down the nonequivalent sign pattern matrices corresponding to each digraph and except for SP \ref{SP4.4}, we use Theorems \ref{thm2}.2 and \ref{thm2}.3 to classify these sign patterns as indicated below.
	\medskip\\
	\begin{minipage}{0.3\textwidth}
		\centering
		\includegraphics[scale=0.6]{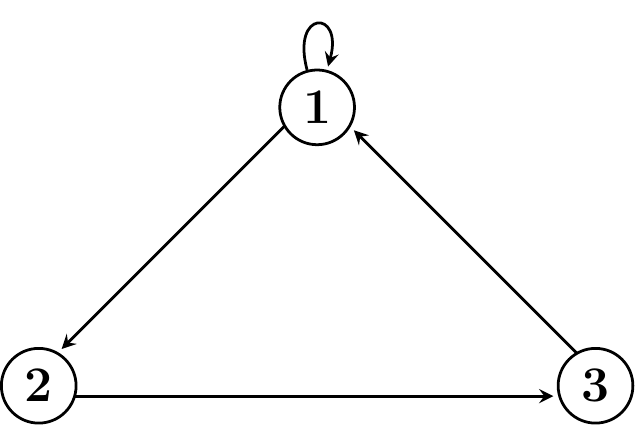}
		\captionof{Digraph}{}
	\end{minipage}
	\begin{minipage}{0.23\textwidth}
		\centering
		$\begin{bmatrix}
		* & + & 0 \\
		0 & 0 & + \\
		+ & 0 & 0
		\end{bmatrix}$
		\captionof{SP}{\textbf{RAP}}
	\end{minipage}
	\begin{minipage}{0.23\textwidth}
		\centering
		$\begin{bmatrix}
		* & - & 0 \\
		0 & 0 & + \\
		+ & 0 & 0
		\end{bmatrix}$
		\captionof{SP}{\textbf{DNA}}
	\end{minipage}
	\begin{minipage}{0.23\textwidth}
		\centering
		$\begin{bmatrix}
		* & + & 0 \\
		0 & 0 & - \\
		+ & 0 & 0
		\end{bmatrix}$
		\captionof{SP}{\textbf{DNA}}
	\end{minipage}
	
	\begin{minipage}{0.3\textwidth}
		\centering
		\includegraphics[scale=0.6]{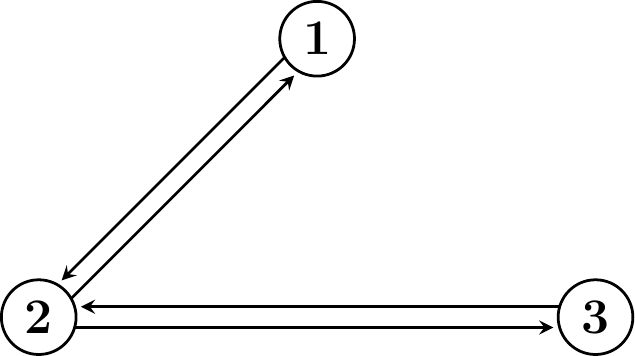}
		\captionof{Digraph}{}
	\end{minipage}
	\begin{minipage}{0.23\textwidth}
		\centering
		$\begin{bmatrix}
		0 & + & 0 \\
		+ & 0 & + \\
		0 & + & 0
		\end{bmatrix}$
		\captionof{SP}{\textbf{RAP}}
	\end{minipage}
	\begin{minipage}{0.23\textwidth}
		\centering
		$\begin{bmatrix}
		0 & - & 0 \\
		+ & 0 & + \\
		0 & + & 0
		\end{bmatrix}$
		\captionof{SP}{\textbf{DNA}}
	\end{minipage}
	\begin{minipage}{0.23\textwidth}
		\centering
		$\begin{bmatrix}
		0 & - & 0 \\
		- & 0 & + \\
		0 & + & 0
		\end{bmatrix}$
		\captionof{SP}{\textbf{DNA}}
		\label{SP3.3}
	\end{minipage}
	\smallskip\\
	\begin{minipage}{0.53\textwidth}
		\centering
		\ 
	\end{minipage}
	\begin{minipage}{0.23\textwidth}
		\centering
		$\begin{bmatrix}
		0 & - & 0 \\
		+ & 0 & - \\
		0 & + & 0
		\end{bmatrix}$
		\captionof{SP}{\textbf{DNA}}
	\end{minipage}
	\begin{minipage}{0.23\textwidth}
		\centering
		$\begin{bmatrix}
		0 & - & 0 \\
		+ & 0 & + \\
		0 & - & 0
		\end{bmatrix}$			
		\captionof{SP}{\textbf{DNA}}
	\end{minipage}
	
	\begin{minipage}{0.3\textwidth}
		\centering
		\includegraphics[scale=0.6]{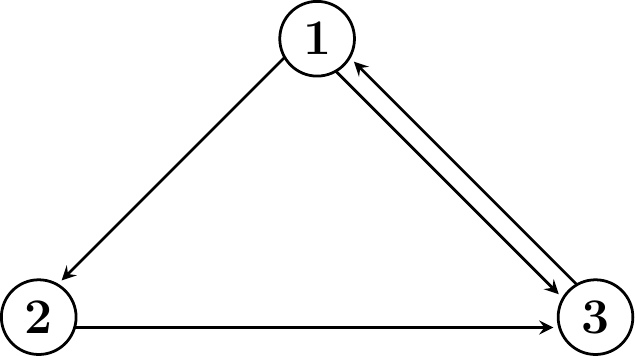}
		\captionof{Digraph}{}
	\end{minipage}
	\begin{minipage}{0.23\textwidth}
		\centering
		$\begin{bmatrix}
		0 & + & + \\
		0 & 0 & + \\
		+ & 0 & 0
		\end{bmatrix}$
		\captionof{SP}{\textbf{RAP}}
	\end{minipage}
	\begin{minipage}{0.23\textwidth}
		\centering
		$\begin{bmatrix}
		0 & + & * \\
		0 & 0 & + \\
		- & 0 & 0
		\end{bmatrix}$
		\captionof{SP}{\textbf{DNA}}
	\end{minipage}
	\begin{minipage}{0.23\textwidth}
		\centering
		$\begin{bmatrix}
		0 & - & * \\
		0 & 0 & + \\
		+ & 0 & 0
		\end{bmatrix}$
		\captionof{SP}{\textbf{DNA}}
	\end{minipage}
	\smallskip\\
	\begin{minipage}{0.3\textwidth}
		\centering
		\ 
	\end{minipage}
	\begin{minipage}{0.23\textwidth}
		\centering
		$\begin{bmatrix}
		0 & + & - \\
		0 & 0 & + \\
		+ & 0 & 0
		\end{bmatrix}$
		\captionof{SP}{\textbf{RAP}}
		\label{SP4.4}
	\end{minipage}
	\begin{minipage}{0.23\textwidth}
		\centering
		$\begin{bmatrix}
		0 & - & + \\
		0 & 0 & - \\
		+ & 0 & 0
		\end{bmatrix}$			
		
		\captionof{SP}{\textbf{DNA}}
	\end{minipage}
	\begin{minipage}{0.23\textwidth}
		\centering
		$\begin{bmatrix}
		0 & - & + \\
		0 & 0 & + \\
		- & 0 & 0
		\end{bmatrix}$			
		
		\captionof{SP}{\textbf{DNA}}
	\end{minipage}\\
	Suppose that the sign pattern of $A=[a_{ij}]$ is SP \ref{SP4.4}. Take $p(x)=k_2x^2+k_1x+k_0$ such that $k_1,k_2>0$, $\frac{k_1}{k_2}> \frac{a_{12}a_{23}}{-a_{13}}$ and $k_0>-k_2a_{13}a_{31}$. Then $p(A)>0$. 
	
	\item \textbf{Graphs with 5 directed edges}. There are six nonequivalent strongly connected digraphs with 5 directed edges. We list down the nonequivalent sign pattern matrices corresponding to each digraph. Except for SP \ref{SP6.4}, SP \ref{SP8.4}, SP \ref{SP9.2}, SP \ref{SP10.4}, and the AAP sign patterns, we use Theorems \ref{thm2}.2,\ref{thm2}.3 and \ref{thm4} to classify these sign patterns as indicated below. 
	\medskip\\
	\begin{minipage}{0.3\textwidth}
		\centering
		\includegraphics[scale=0.6]{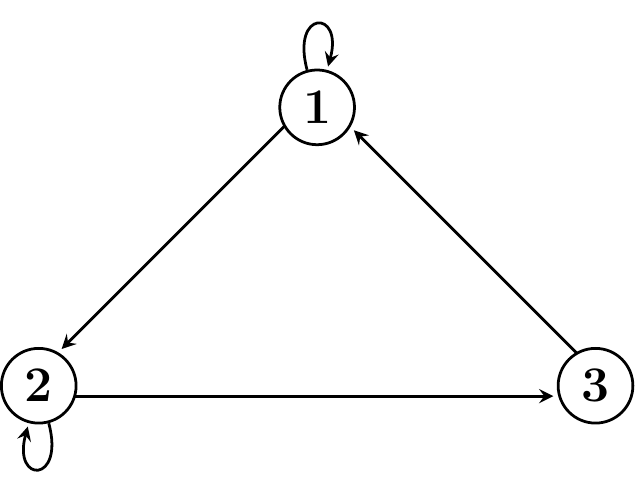}
		\captionof{Digraph}{}
	\end{minipage}
	\begin{minipage}{0.23\textwidth}
		\centering
		$\begin{bmatrix}
		* & + & 0 \\
		0 & * & + \\
		+ & 0 & 0
		\end{bmatrix}$
		\captionof{SP}{\textbf{RAP}}
	\end{minipage}
	\begin{minipage}{0.23\textwidth}
		\centering
		$\begin{bmatrix}
		* & - & 0 \\
		0 & * & + \\
		+ & 0 & 0
		\end{bmatrix}$
		\captionof{SP}{\textbf{DNA}}
	\end{minipage}
	\begin{minipage}{0.23\textwidth}
		\centering
		$\begin{bmatrix}
		* & + & 0 \\
		0 & * & - \\
		+ & 0 & 0
		\end{bmatrix}$
		\captionof{SP}{\textbf{DNA}}
	\end{minipage}
	
	\begin{minipage}{0.3\textwidth}
		\centering
		\includegraphics[scale=0.6]{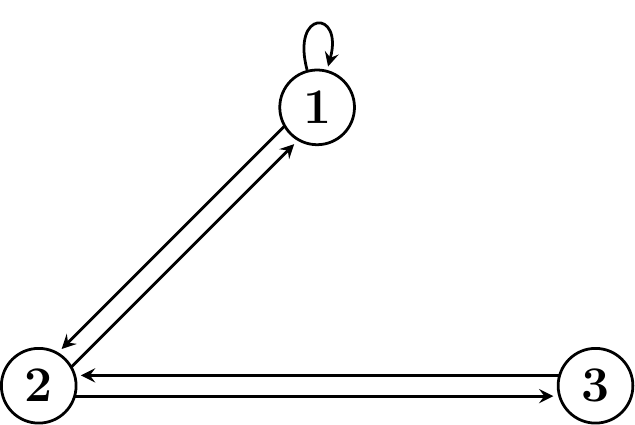}
		\captionof{Digraph}{}
	\end{minipage}
	\begin{minipage}{0.23\textwidth}
		\centering
		$\begin{bmatrix}
		* & + & 0 \\
		+ & 0 & + \\
		0 & + & 0
		\end{bmatrix}$
		\captionof{SP}{\textbf{RAP}}
	\end{minipage}
	\begin{minipage}{0.23\textwidth}
		\centering
		$\begin{bmatrix}
		* & - & 0 \\
		+ & 0 & + \\
		0 & + & 0
		\end{bmatrix}$
		\captionof{SP}{\textbf{DNA}}
	\end{minipage}
	\begin{minipage}{0.23\textwidth}
		\centering
		$\begin{bmatrix}
		* & + & 0 \\
		+ & 0 & - \\
		0 & + & 0
		\end{bmatrix}$
		\captionof{SP}{\textbf{DNA}}
	\end{minipage}
	\smallskip\\
	\begin{minipage}{0.3\textwidth}
		\centering
		\ 
	\end{minipage}
	\begin{minipage}{0.23\textwidth}
		\centering
		$\begin{bmatrix}
		+ & - & 0 \\
		- & 0 & + \\
		0 & + & 0
		\end{bmatrix}$
		\captionof{SP}{\textbf{RAP}}
		\label{SP6.4}
	\end{minipage}
	\begin{minipage}{0.23\textwidth}
		\centering
		$\begin{bmatrix}
		+ & + & 0 \\
		+ & 0 & - \\
		0 & - & 0
		\end{bmatrix}$			
		
		\captionof{SP}{\textbf{DNA}}
	\end{minipage}
	\begin{minipage}{0.23\textwidth}
		\centering
		$\begin{bmatrix}
		* & - & 0 \\
		+ & 0 & - \\
		0 & + & 0
		\end{bmatrix}$			
		
		\captionof{SP}{\textbf{DNA}}
	\end{minipage}\medskip\\
Suppose that the sign pattern of $A=[a_{ij}]$ is SP \ref{SP6.4}. Take $p(x)=k_2x^2+k_1x+k_0$ such that $k_1>0>k_2$ and $\frac{k_1}{k_2}> -a_{11}$, $k_0>\max\{-k_2(a_{11}^2+a_{12}a_{21})-k_1a_{11},-k_2(a_{12}a_{21}+a_{31}a_{23}), -k_2a_{31}a_{23}\}$. Then $p(A)>0$. \medskip\\
	\begin{minipage}{0.3\textwidth}
		\centering
		\includegraphics[scale=0.6]{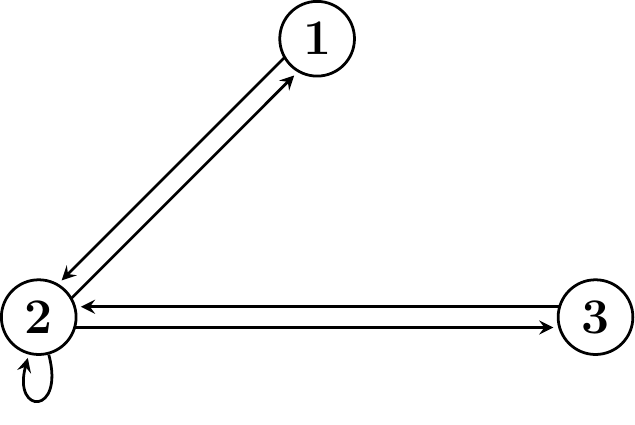}
		\captionof{Digraph}{}
	\end{minipage}
	\begin{minipage}{0.23\textwidth}
		\centering
		$\begin{bmatrix}
		0 & + & 0 \\
		+ & * & + \\
		0 & + & 0
		\end{bmatrix}$
		\captionof{SP}{\textbf{RAP}}
		\label{SP7.1}
	\end{minipage}
	\begin{minipage}{0.23\textwidth}
		\centering
		$\begin{bmatrix}
		0 & - & 0 \\
		+ & * & + \\
		0 & + & 0
		\end{bmatrix}$
		\captionof{SP}{\textbf{DNA}}
	\end{minipage}
	\begin{minipage}{0.23\textwidth}
		\centering
		$\begin{bmatrix}
		0 & - & 0 \\
		- & * & + \\
		0 & + & 0
		\end{bmatrix}$
		\captionof{SP}{\textbf{DNA}}
	\end{minipage}
	\smallskip\\
	\begin{minipage}{0.53\textwidth}
		\centering
		\ 
	\end{minipage}
	\begin{minipage}{0.23\textwidth}
		\centering
		$\begin{bmatrix}
		0 & - & 0 \\
		+ & * & - \\
		0 & + & 0
		\end{bmatrix}$
		\captionof{SP}{\textbf{DNA}}
	\end{minipage}
	\begin{minipage}{0.23\textwidth}
		\centering
		$\begin{bmatrix}
		0 & - & 0 \\
		+ & * & + \\
		0 & - & 0
		\end{bmatrix}$			
		\captionof{SP}{\textbf{DNA}}
	\end{minipage}
	
	\begin{minipage}{0.3\textwidth}
		\centering
		\includegraphics[scale=0.6]{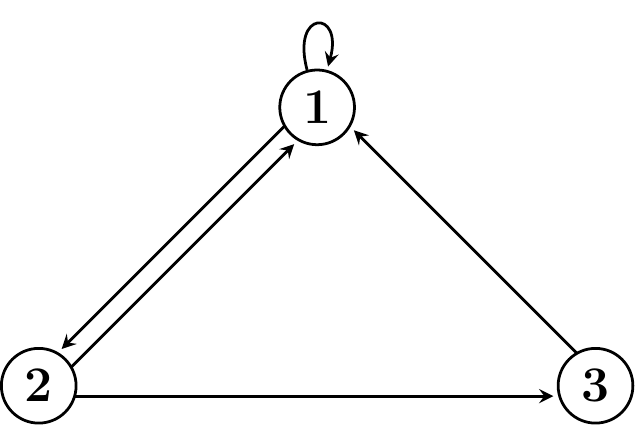}
		\captionof{Digraph}{}
		\label{D8}
	\end{minipage}
	\begin{minipage}{0.23\textwidth}
		\centering
		$\begin{bmatrix}
		* & + & 0 \\
		+ & 0 & + \\
		+ & 0 & 0
		\end{bmatrix}$
		\captionof{SP}{\textbf{RAP}}
	\end{minipage}
	\begin{minipage}{0.23\textwidth}
		\centering
		$\begin{bmatrix}
		+ & + & 0 \\
		- & 0 & + \\
		+ & 0 & 0
		\end{bmatrix}$
		\captionof{SP}{\textbf{AAP}}
		\label{SP8.2}
	\end{minipage}
	\begin{minipage}{0.23\textwidth}
		\centering
		$\begin{bmatrix}
		* & - & 0 \\
		* & 0 & + \\
		+ & 0 & 0
		\end{bmatrix}$
		\captionof{SP}{\textbf{DNA}}
	\end{minipage}
	\smallskip\\
	\begin{minipage}{0.3\textwidth}
		\centering
		\ 
	\end{minipage}
	\begin{minipage}{0.23\textwidth}
		\centering
		$\begin{bmatrix}
		- & + & 0 \\
		- & 0 & + \\
		+ & 0 & 0
		\end{bmatrix}$
		\captionof{SP}{\textbf{RAP}}
		\label{SP8.4}
	\end{minipage}
	\begin{minipage}{0.23\textwidth}
		\centering
		$\begin{bmatrix}
		* & + & 0 \\
		* & 0 & - \\
		+ & 0 & 0
		\end{bmatrix}$				
		\captionof{SP}{\textbf{DNA}}
	\end{minipage}
	\begin{minipage}{0.23\textwidth}
		\centering
		$\begin{bmatrix}
		* & + & 0 \\
		* & 0 & + \\
		- & 0 & 0
		\end{bmatrix}$			
		\captionof{SP}{\textbf{DNA}}
	\end{minipage}
	
Suppose that the digraph of $A=[a_{ij}]\in M_3(\mathbb{R})$ is Digraph \ref{D8}. If the sign pattern of $A$ is
\begin{itemize}
\item  SP \ref{SP8.2}, then $A$ is AP if and only if $a_{11}a_{21}+a_{23}a_{31}>0$. In particular, if we take $a_{11}=a_{12}=a_{23}=a_{31}=-a_{21}=1$, then $A$ is not AP. On the other hand, if we take $a_{11}=a_{12}=-a_{21}=1$ and $a_{23}=a_{31}=10$, then $A$ is AP.
\item SP \ref{SP8.4}, take $p(x)=k_2x^2+k_1x+k_0$ such that $k_1,k_2>0$, $\frac{k_1}{k_2}>-a_{11}-\dfrac{a_{23}a_{31}}{a_{21}}$ and $k_0>-k_1a_{11}-k_2(a_{11}^2+a_{12}a_{21})$ so that $p(A)>0$.
\end{itemize}

	\begin{minipage}{0.3\textwidth}
		\centering
		\includegraphics[scale=0.6]{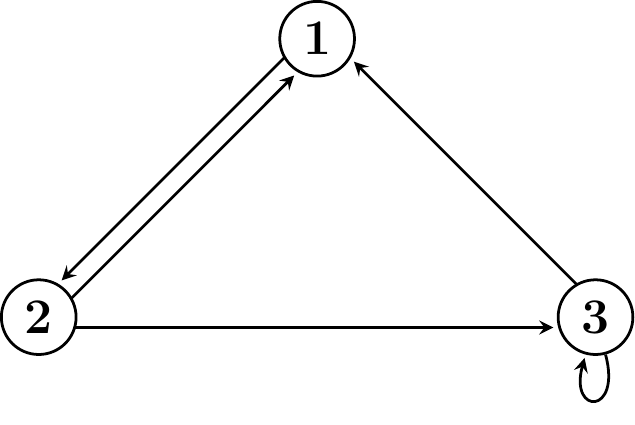}
		\captionof{Digraph}{}
		\label{D9}
	\end{minipage}
	\begin{minipage}{0.23\textwidth}
		\centering
		$\begin{bmatrix}
		0 & + & 0 \\
		+ & 0 & + \\
		+ & 0 & *
		\end{bmatrix}$
		\captionof{SP}{\textbf{RAP}}
	\end{minipage}
	\begin{minipage}{0.23\textwidth}
		\centering
		$\begin{bmatrix}
		0 & + & 0 \\
		- & 0 & + \\
		+ & 0 & +
		\end{bmatrix}$
		\captionof{SP}{\textbf{RAP}}
		\label{SP9.2}
	\end{minipage}
	\begin{minipage}{0.23\textwidth}
		\centering
		$\begin{bmatrix}
		0 & + & 0 \\
		- & 0 & + \\
		+ & 0 & -
		\end{bmatrix}$
		\captionof{SP}{\textbf{AAP}}
		\label{SP9.3}
	\end{minipage}
	\smallskip\\
	\begin{minipage}{0.3\textwidth}
		\centering
		\ 
	\end{minipage}
	\begin{minipage}{0.23\textwidth}
		\centering
		$\begin{bmatrix}
		0 & - & 0 \\
		+ & 0 & + \\
		+ & 0 & *
		\end{bmatrix}$
		\captionof{SP}{\textbf{DNA}}
	\end{minipage}
	\begin{minipage}{0.23\textwidth}
		\centering
		$\begin{bmatrix}
		0 & + & 0 \\
		+ & 0 & - \\
		+ & 0 & *
		\end{bmatrix}$				
		\captionof{SP}{\textbf{DNA}}
	\end{minipage}
	\begin{minipage}{0.23\textwidth}
		\centering
		$\begin{bmatrix}
		0 & + & 0 \\
		+ & 0 & - \\
		- & 0 & +
		\end{bmatrix}$			
		\captionof{SP}{\textbf{AAP}}
		\label{SP9.6}
	\end{minipage}
	\smallskip\\
	\begin{minipage}{0.3\textwidth}
		\centering
		\ 
	\end{minipage}
	\begin{minipage}{0.23\textwidth}
		\centering
		$\begin{bmatrix}
		0 & - & 0 \\
		+ & 0 & - \\
		+ & 0 & *
		\end{bmatrix}$
		\captionof{SP}{\textbf{DNA}}
	\end{minipage}
	\begin{minipage}{0.23\textwidth}
		\centering
		$\begin{bmatrix}
		0 & - & 0 \\
		- & 0 & + \\
		+ & 0 & +
		\end{bmatrix}$				
		\captionof{SP}{\textbf{DNA}}
	\end{minipage}
	\begin{minipage}{0.23\textwidth}
		\centering
		$\begin{bmatrix}
		0 & + & 0 \\
		- & 0 & - \\
		+ & 0 & *
		\end{bmatrix}$			
		\captionof{SP}{\textbf{DNA}}
	\end{minipage}
	
	Suppose that the digraph of $A=[a_{ij}]\in M_3(\mathbb{R})$ is Digraph \ref{D9}. If the sign pattern of $A$ is
	\begin{itemize}
	\item SP \ref{SP9.2}, take $p(x)=k_2x^2+k_1x+k_0$ such that $k_1,k_2>0$, $\frac{k_1}{k_2}>\frac{a_{31}a_{23}}{-a_{21}}$ and $k_0>-k_2a_{21}a_{12}$ so that $p(A)>0$.
	\item SP \ref{SP9.3}, then $A$ is AP if and only if $a_{21}a_{33} < a_{23}a_{31}$. In particular, if we take $a_{12}=a_{23}=a_{31}=-a_{21}=-a_{33}=1$, then $A$ is not AP. On the other hand, if we take $a_{12}=a_{23}=a_{31}=2$ and $a_{21}=a_{33}=-1$, then $A$ is AP.
		\item SP \ref{SP9.6}, then $A$ is AP if and only if $a_{21}a_{33} > a_{23}a_{31}$. In particular, if we take $a_{12}=-a_{23}=-a_{31}=a_{21}=a_{33}=1$, then $A$ is not AP. On the other hand, if we take $a_{12}=a_{21}=a_{33}=2$ and $a_{23}=a_{31}=-1$, then $A$ is AP.
	\end{itemize} 
	
	\begin{minipage}{0.3\textwidth}
		\centering
		\includegraphics[scale=0.6]{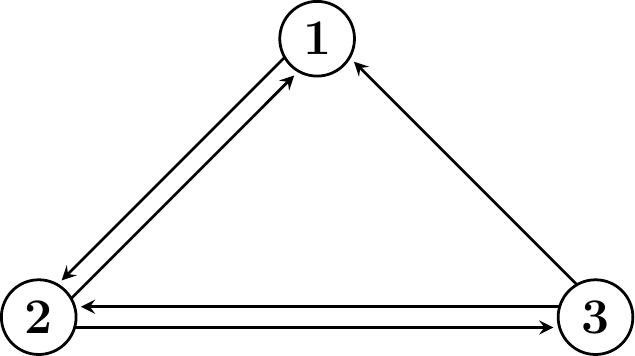}
		\captionof{Digraph}{}
		\label{D10}
	\end{minipage}
	\begin{minipage}{0.23\textwidth}
		\centering
		$\begin{bmatrix}
		0 & + & 0 \\
		+ & 0 & + \\
		+ & + & 0
		\end{bmatrix}$
		\captionof{SP}{\textbf{RAP}}
	\end{minipage}
	\begin{minipage}{0.23\textwidth}
		\centering
		$\begin{bmatrix}
		0 & + & 0 \\
		+ & 0 & + \\
		+ & - & 0
		\end{bmatrix}$
		\captionof{SP}{\textbf{AAP}}
		\label{SP10.2}
	\end{minipage}
	\begin{minipage}{0.23\textwidth}
		\centering
		$\begin{bmatrix}
		0 & + & 0 \\
		* & 0 & - \\
		* & * & 0
		\end{bmatrix}$
		\captionof{SP}{\textbf{DNA}}
	\end{minipage}
	\smallskip\\
	\begin{minipage}{0.3\textwidth}
		\centering
		\ 
	\end{minipage}
	\begin{minipage}{0.23\textwidth}
		\centering
		$\begin{bmatrix}
		0 & + & 0 \\
		- & 0 & + \\
		+ & - & 0
		\end{bmatrix}$
		\captionof{SP}{\textbf{RAP}}
		\label{SP10.4}
	\end{minipage}
	\begin{minipage}{0.23\textwidth}
		\centering
		$\begin{bmatrix}
		0 & + & 0 \\
		+ & 0 & + \\
		- & + & 0
		\end{bmatrix}$				
		\captionof{SP}{\textbf{AAP}}
		\label{SP10.5}
	\end{minipage}
	\begin{minipage}{0.23\textwidth}
		\centering
		$\begin{bmatrix}
		0 & + & 0 \\
		- & 0 & + \\
		- & + & 0
		\end{bmatrix}$			
		\captionof{SP}{\textbf{DNA}}
	\end{minipage}
	\smallskip\\
	\begin{minipage}{0.3\textwidth}
		\centering
		\ 
	\end{minipage}
	\begin{minipage}{0.23\textwidth}
		\centering
		\ 
	\end{minipage}
	\begin{minipage}{0.23\textwidth}
		\centering
		\ 
	\end{minipage}
	\begin{minipage}{0.23\textwidth}
		\centering
		$\begin{bmatrix}
		0 & - & 0 \\
		+ & 0 & - \\
		+ & + & 0
		\end{bmatrix}$			
		\captionof{SP}{\textbf{DNA}}
	\end{minipage}
	
	Suppose that the digraph of $A=[a_{ij}]\in M_3(\mathbb{R})$ is Digraph \ref{D10}. If the sign pattern of $A$ is
	\begin{itemize}
		\item SP \ref{SP10.2}, then $A$ is AP if and only if $a_{21}a_{32}^2 < a_{12}a_{31}^2$. In particular, if we take $a_{12}=a_{21}=a_{23}=a_{31}=-a_{32}=1$, then $A$ is not AP. On the other hand, if we take $a_{12}=a_{21}=a_{23}=a_{31}=2$ and $a_{32}=-1$, then $A$ is AP.
		\item SP \ref{SP10.4}, take $p(x)=k_2x^2+k_1x+k_0$ such that $k_1,k_2>0$, $\frac{k_1}{k_2}>\max\left\{\frac{a_{31}a_{12}}{-a_{32}}, \frac{a_{31}a_{23}}{-a_{21}}\right\}$ and $k_0> -k_2(a_{12}a_{21}+a_{23}a_{32})$ so that $p(A)>0$. 
		\item SP \ref{SP10.5}, then $A$ is AP if and only if $a_{23}a_{31}^2 < a_{21}^2a_{32}$ and $a_{21}a_{31}^2 < a_{21}a_{32}^2$. In particular, if we take $a_{12}=a_{21}=a_{23}=-a_{31}=a_{32}=1$, then $A$ is not AP. On the other hand, if we take $a_{12}=a_{21}=a_{23}=a_{32}=2$ and $a_{31}=-1$, then $A$ is AP.
	\end{itemize}

	\item \textbf{Graphs with 6 directed edges}. There are eight nonequivalent strongly connected digraphs with 6 directed edges. We list down the nonequivalent sign pattern matrices corresponding to each digraph and except for SP 12.4, SP \ref{SP13.4}, SP \ref{SP14.4}, SP \ref{SP15.4}, SP \ref{SP17.4}, SP \ref{SP18.4}, and the AAP sign patterns, we use Theorems \ref{thm2}.2,\ref{thm2}.3 and \ref{thm4} to classify these sign patterns as indicated below. 
	
	\begin{minipage}{0.4\textwidth}
		\centering
		\includegraphics[scale=0.6]{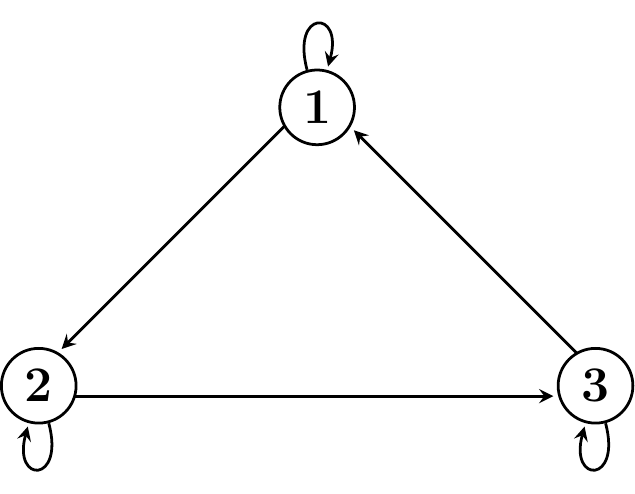}
		\captionof{Digraph}{}
	\end{minipage}
	\begin{minipage}{0.28\textwidth}
		\centering 
		$\begin{bmatrix} 
		* & + & 0 \\ 
		0 & * & + \\
		+ & 0 & * 
		\end{bmatrix}$
		\captionof{SP}{\textbf{RAP}}
	\end{minipage}
	\begin{minipage}{0.28\textwidth}
		\centering
		$\begin{bmatrix} 
		* & - & 0 \\
		0 & * & + \\ 
		+ & 0 & * 
		\end{bmatrix}$
		\captionof{SP}{\textbf{DNA}}
	\end{minipage}
	
	\begin{minipage}{0.3\textwidth}
		\centering
		\includegraphics[scale=0.6]{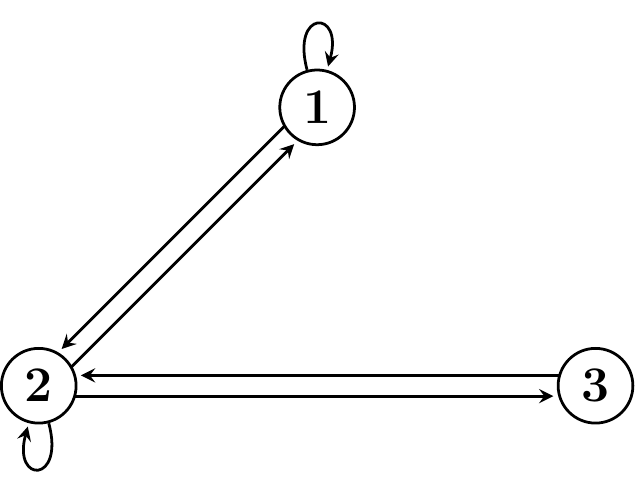}
		\captionof{Digraph}{}
	\end{minipage}
	\begin{minipage}{0.23\textwidth}
		\centering
		$\begin{bmatrix}
		* & + & 0 \\
		+ & * & + \\
		0 & + & 0
		\end{bmatrix}$
		\captionof{SP}{\textbf{RAP}}
	\end{minipage}
	\begin{minipage}{0.23\textwidth}
		\centering
		$\begin{bmatrix}
		* & + & 0 \\
		- & * & + \\
		0 & + & 0
		\end{bmatrix}$
		\captionof{SP}{\textbf{DNA}}
	\end{minipage}
	\begin{minipage}{0.23\textwidth}
		\centering
		$\begin{bmatrix}
		* & + & 0 \\
		* & * & + \\
		0 & - & 0
		\end{bmatrix}$
		\captionof{SP}{\textbf{DNA}}
	\end{minipage}
	\smallskip\\
	\begin{minipage}{0.53\textwidth}
		\centering
		\ 
	\end{minipage}
	\begin{minipage}{0.23\textwidth}
		\centering
		$\begin{bmatrix}
		+ & - & 0 \\
		- & * & + \\
		0 & + & 0
		\end{bmatrix}$
		\captionof{SP}{\textbf{RAP}}
		\label{SP12.4}
	\end{minipage}
	\begin{minipage}{0.23\textwidth}
		\centering
		$\begin{bmatrix}
		+ & + & 0 \\
		+ & * & - \\
		0 & - & 0
		\end{bmatrix}$			
		\captionof{SP}{\textbf{DNA}}
	\end{minipage}\medskip\\
 Suppose that the sign pattern of $A=[a_{ij}]$ is  \ref{SP12.4}. Then take $p(x)=k_2x^2+k_1x+k_0$ such that $k_2<0$, $-(a_{11}+a_{22})< \frac{k_1}{k_2}<-a_{22}$ and $k_0$ be larger than all the diagonal entries of $-k_1A-k_2A^2$ so that $p(A)>0$. 
  \medskip\\
	\begin{minipage}{0.3\textwidth}
		\centering
		\includegraphics[scale=0.6]{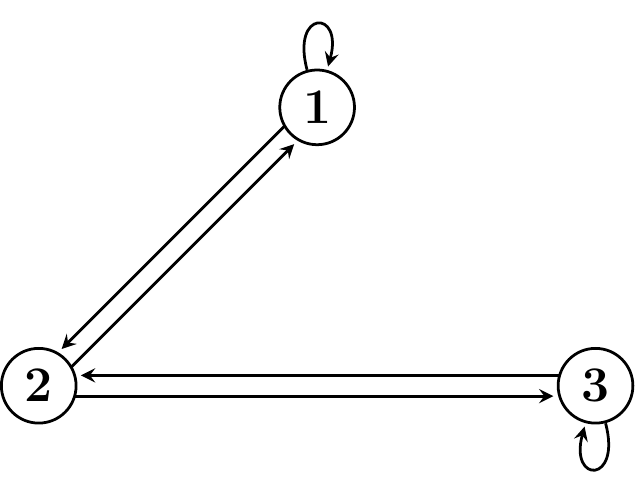}
		\captionof{Digraph}{}
		\label{D13}
	\end{minipage}
	\begin{minipage}{0.23\textwidth}
		\centering
		$\begin{bmatrix}
		* & + & 0 \\
		+ & 0 & + \\
		0 & + & *
		\end{bmatrix}$
		\captionof{SP}{\textbf{RAP}}
	\end{minipage}
	\begin{minipage}{0.23\textwidth}
		\centering
		$\begin{bmatrix}
		+ & + & 0 \\
		+ & 0 & - \\
		0 & - & +
		\end{bmatrix}$
		\captionof{SP}{\textbf{AAP}}
		\label{SP13.2}
	\end{minipage}
	\begin{minipage}{0.23\textwidth}
		\centering
		$\begin{bmatrix}
		* & + & 0 \\
		+ & 0 & + \\
		0 & - & *
		\end{bmatrix}$
		\captionof{SP}{\textbf{DNA}}
	\end{minipage}
	\smallskip\\
	\begin{minipage}{0.3\textwidth}
		\centering
		\ 
	\end{minipage}
	\begin{minipage}{0.23\textwidth}
		\centering
		$\begin{bmatrix}
		+ & - & 0 \\
		- & 0 & + \\
		0 & + & -
		\end{bmatrix}$
		\captionof{SP}{\textbf{RAP}}
		\label{SP13.4}
	\end{minipage}
	\begin{minipage}{0.23\textwidth}
		\centering
		$\begin{bmatrix}
		* & - & 0 \\
		+ & 0 & + \\
		0 & - & *
		\end{bmatrix}$				
		\captionof{SP}{\textbf{DNA}}
	\end{minipage}
	\begin{minipage}{0.23\textwidth}
		\centering
		$\begin{bmatrix}
		- & - & 0 \\
		- & 0 & + \\
		0 & + & +
		\end{bmatrix}$			
		\captionof{SP}{\textbf{DNA}}
		\label{SP13.6}
	\end{minipage}
	
	Suppose that the digraph of $A=[a_{ij}]\in M_3(\mathbb{R})$ is Digraph \ref{D13}. If the sign pattern of $A$ is 
\begin{itemize}
\item SP \ref{SP13.2}, then $A$ is AP if and only if $a_{33}>a_{11}$. In particular, if we take $a_{11}=a_{12}=a_{21}=-a_{23}=-a_{32}=a_{33}=1$, then $A$ is not AP. On the other hand, if we take $a_{11}=a_{12}=a_{21}=-a_{23}=-a_{32}=1$ and $a_{33}=3$, then $A$ is AP.
\item SP \ref{SP13.4}, then take $p(x)=k_2x^2+k_1x+k_0$ such that $k_2<0$, $-a_{11}< \frac{k_1}{k_2}<-a_{33}$ and $k_0$ be larger than all the diagonal entries of $-k_1A-k_2A^2$ so that $p(A)>0$. 
\end{itemize}	
	
	\begin{minipage}{0.3\textwidth}
		\centering
		\includegraphics[scale=0.6]{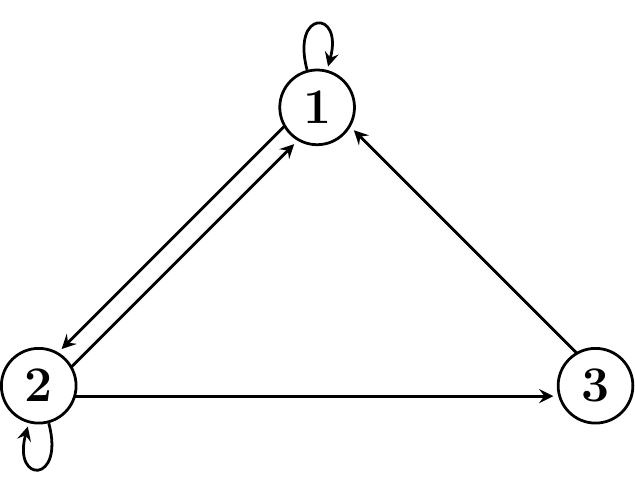}
		\captionof{Digraph}{}
		\label{D14}
	\end{minipage}
	\begin{minipage}{0.23\textwidth}
		\centering
		$\begin{bmatrix}
		* & + & 0 \\
		+ & * & + \\
		+ & 0 & 0
		\end{bmatrix}$
		\captionof{SP}{\textbf{RAP}}
	\end{minipage}
	\begin{minipage}{0.23\textwidth}
		\centering
		$\begin{bmatrix}
		+ & + & 0 \\
		- & + & + \\
		+ & 0 & 0
		\end{bmatrix}$
		\captionof{SP}{\textbf{AAP}}
		\label{SP14.2}
	\end{minipage}
	\begin{minipage}{0.23\textwidth}
		\centering
		$\begin{bmatrix}
		- & + & 0 \\
		- & + & + \\
		+ & 0 & 0
		\end{bmatrix}$
		\captionof{SP}{\textbf{AAP}}
		\label{SP14.3}
	\end{minipage}
	\smallskip\\
	\begin{minipage}{0.3\textwidth}
		\centering
		\ 
	\end{minipage}
	\begin{minipage}{0.23\textwidth}
		\centering
		$\begin{bmatrix}
		- & + & 0 \\
		- & - & + \\
		+ & 0 & 0
		\end{bmatrix}$
		\captionof{SP}{\textbf{RAP}}
			\label{SP14.4}
	\end{minipage}
	\begin{minipage}{0.23\textwidth}
		\centering
		$\begin{bmatrix}
		+ & - & 0 \\
		- & + & + \\
		+ & 0 & 0
		\end{bmatrix}$				
		\captionof{SP}{\textbf{AAP}}
		\label{SP14.5}
	\end{minipage}
	\begin{minipage}{0.23\textwidth}
		\centering
		$\begin{bmatrix}
		* & - & 0 \\
		+ & * & + \\
		+ & 0 & 0
		\end{bmatrix}$			
		\captionof{SP}{\textbf{DNA}}
	\end{minipage}
	\smallskip\\
	\begin{minipage}{0.3\textwidth}
		\centering
		\ 
	\end{minipage}
	\begin{minipage}{0.23\textwidth}
		\centering
		$\begin{bmatrix}
		* & + & 0 \\
		+ & * & - \\
		+ & 0 & 0
		\end{bmatrix}$
		\captionof{SP}{\textbf{DNA}}
	\end{minipage}
	\begin{minipage}{0.23\textwidth}
		\centering
		$\begin{bmatrix}
		* & - & 0 \\
		+ & * & + \\
		- & 0 & 0
		\end{bmatrix}$				
		\captionof{SP}{\textbf{DNA}}
	\end{minipage}
	\begin{minipage}{0.23\textwidth}
		\centering
		$\begin{bmatrix}
		* & + & 0 \\
		+ & * & - \\
		- & 0 & 0
		\end{bmatrix}$			
		\captionof{SP}{\textbf{DNA}}
	\end{minipage}
	
	Suppose that the digraph of $A=[a_{ij}]\in M_3(\mathbb{R})$ is Digraph \ref{D14}. If the sign pattern of $A$ is
	\begin{itemize}
		\item SP \ref{SP14.2}, then $A$ is AP if and only if $a_{11}a_{21}, a_{22}a_{21} > -a_{23}a_{31}$. In particular, if we take $a_{11}=a_{12}=-a_{21}=a_{22}=a_{23}=a_{31}=1$, then $A$ is not AP. On the other hand, if we take $a_{11}=a_{12}=a_{22}=a_{23}=a_{31}=2$ and $a_{21}=-1$, then $A$ is AP.
		\item SP \ref{SP14.3}, then $A$ is AP if and only if $a_{22}a_{21} > -a_{23}a_{31}$. In particular, if we take $-a_{11}=a_{12}=-a_{21}=a_{22}=a_{23}=a_{31}=1$, then $A$ is not AP. On the other hand, if we take $a_{12}=a_{22}=a_{23}=a_{31}=2$ and $a_{11}=a_{21}=-1$, then $A$ is AP.
		\item SP \ref{SP14.4}, then take $p(x)=k_2x^2+k_1x+k_0$ such that $k_1,k_2>0$, $-(a_{11}+a_{22})< \frac{k_1}{k_2}<-(a_{11}+a_{22})-\frac{a_{23}a_{31}}{a_{21}}$ and $k_0$ be larger than all the diagonal entries of $-k_1A-k_2A^2$ so that $p(A)>0$. 
		\item SP \ref{SP14.5}, then $A$ is AP if and only if $a_{11}a_{21}, a_{22}a_{21} < -a_{23}a_{31}$. In particular, if we take $a_{11}=-a_{12}=-a_{21}=a_{22}=a_{23}=a_{31}=1$, then $A$ is not AP. On the other hand, if we take $a_{11}=-a_{12}=a_{22}=a_{23}=a_{31}=1$ and $a_{21}=-2$, then $A$ is AP.
	\end{itemize} 
	
	\begin{minipage}{0.3\textwidth}
		\centering
		\includegraphics[scale=0.6]{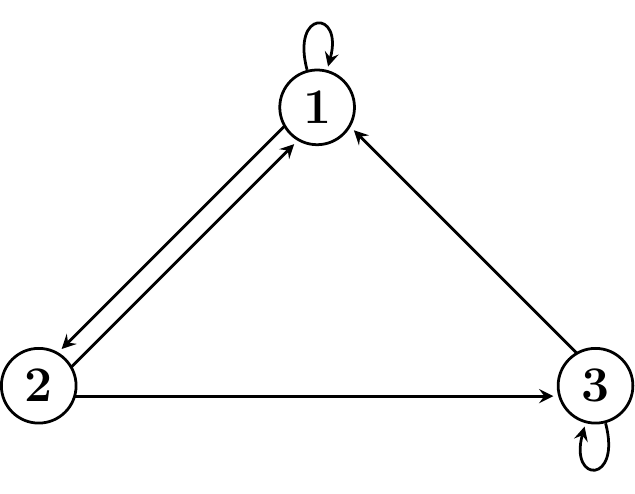}
		\captionof{Digraph}{}
		\label{D15}
	\end{minipage}
	\begin{minipage}{0.23\textwidth}
		\centering
		$\begin{bmatrix}
		* & + & 0 \\
		+ & 0 & + \\
		+ & 0 & *
		\end{bmatrix}$
		\captionof{SP}{\textbf{RAP}}
	\end{minipage}
	\begin{minipage}{0.23\textwidth}
		\centering
		$\begin{bmatrix}
		+ & + & 0 \\
		- & 0 & + \\
		+ & 0 & +
		\end{bmatrix}$
		\captionof{SP}{\textbf{AAP}}
		\label{SP15.2}
	\end{minipage}
	\begin{minipage}{0.23\textwidth}
		\centering
		$\begin{bmatrix}
		* & + & 0 \\
		* & 0 & - \\
		+ & 0 & *
		\end{bmatrix}$
		\captionof{SP}{\textbf{DNA}}
	\end{minipage}
	\smallskip\\
	\begin{minipage}{0.3\textwidth}
		\centering
		\ 
	\end{minipage}
	\begin{minipage}{0.23\textwidth}
		\centering
		$\begin{bmatrix}
		- & + & 0 \\
		- & 0 & + \\
		+ & 0 & +
		\end{bmatrix}$
		\captionof{SP}{\textbf{RAP}}
		\label{SP15.4}
	\end{minipage}
	\begin{minipage}{0.23\textwidth}
		\centering
		$\begin{bmatrix}
		+ & + & 0 \\
		- & 0 & + \\
		+ & 0 & -
		\end{bmatrix}$				
		\captionof{SP}{\textbf{AAP}}
		\label{SP15.5}
	\end{minipage}
	\begin{minipage}{0.23\textwidth}
		\centering
		$\begin{bmatrix}
		* & - & 0 \\
		+ & 0 & + \\
		+ & 0 & *
		\end{bmatrix}$			
		\captionof{SP}{\textbf{DNA}}
	\end{minipage}
	\smallskip\\
	\begin{minipage}{0.3\textwidth}
		\centering
		\ 
	\end{minipage}
	\begin{minipage}{0.23\textwidth}
		\centering
		$\begin{bmatrix}
		- & + & 0 \\
		- & 0 & + \\
		+ & 0 & -
		\end{bmatrix}$
		\captionof{SP}{\textbf{AAP}}
		\label{SP15.7}
	\end{minipage}
	\begin{minipage}{0.23\textwidth}
		\centering
		$\begin{bmatrix}
		+ & + & 0 \\
		+ & 0 & - \\
		- & 0 & +
		\end{bmatrix}$				
		\captionof{SP}{\textbf{AAP}}
		\label{SP15.8}
	\end{minipage}
	\begin{minipage}{0.23\textwidth}
		\centering
		$\begin{bmatrix}
		* & + & 0 \\
		* & 0 & + \\
		- & 0 & *
		\end{bmatrix}$			
		\captionof{SP}{\textbf{DNA}}
	\end{minipage}
	\smallskip\\
	\begin{minipage}{0.53\textwidth}
		\centering
		\ 
	\end{minipage}
	\begin{minipage}{0.23\textwidth}
		\centering
		$\begin{bmatrix}
		- & + & 0 \\
		+ & 0 & - \\
		- & 0 & +
		\end{bmatrix}$
		\captionof{SP}{\textbf{AAP}}
		\label{SP15.10}
	\end{minipage}
	\begin{minipage}{0.23\textwidth}
		\centering
		$\begin{bmatrix}
		* & + & 0 \\
		+ & 0 & - \\
		- & 0 & -
		\end{bmatrix}$			
		\captionof{SP}{\textbf{DNA}}
	\end{minipage}\\
	
	Suppose that the digraph of $A=[a_{ij}]\in M_3(\mathbb{R})$ is Digraph \ref{D15}. If the sign pattern of $A$ is
	\begin{itemize}
		\item SP \ref{SP15.2}, then $A$ is AP if and only if $a_{11}a_{21}- a_{33}a_{21} > -a_{23}a_{31}$. In particular, if we take $a_{12}=-a_{21}=a_{23}=a_{31}=a_{33}=1$ and $a_{11}=2$, then $A$ is not AP. On the other hand, if we take $a_{11}=a_{12}=-a_{21}=a_{23}=a_{31}=a_{33}=1$, then $A$ is AP.
		\item SP \ref{SP15.4}, then take $p(x)=k_2x^2+k_1x+k_0$ such that $k_1,k_2>0$, $-a_{11}< \frac{k_1}{k_2}<-a_{11}-\frac{a_{23}a_{31}}{a_{21}}$ and $k_0$ be larger than all the diagonal entries of $-k_1A-k_2A^2$ so that $p(A)>0$. 
		\item SP \ref{SP15.5}, then $A$ is AP if and only if $a_{11}a_{21}- a_{33}a_{21} > -a_{23}a_{31}$. In particular, if we take $a_{11}=a_{12}=-a_{21}=a_{23}=a_{31}=-a_{33}=1$, then $A$ is not AP. On the other hand, if we take $a_{11}=a_{12}=a_{23}=a_{31}=2$ and $a_{21}=a_{33}=-1$, then $A$ is AP.
		\item SP \ref{SP15.7}, then $A$ is AP if and only if $a_{33}a_{21} < a_{23}a_{31}$. In particular, if we take $-a_{11}=a_{12}=-a_{21}=a_{23}=a_{31}=-a_{33}=1$, then $A$ is not AP. On the other hand, if we take $a_{12}=a_{23}=a_{31}=2$ and $a_{11}=a_{21}=a_{33}=-1$, then $A$ is AP.
		\item SP \ref{SP15.8}, then $A$ is AP if and only if $a_{11}a_{21}- a_{33}a_{21} < -a_{23}a_{31}$. In particular, if we take $a_{11}=a_{12}=a_{21}=-a_{23}=-a_{31}=a_{33}=1$, then $A$ is not AP. On the other hand, if we take $a_{11}=a_{12}=a_{21}=-a_{23}=-a_{31}=1$ and $a_{33}=3$, then $A$ is AP.
		\item SP \ref{SP15.10}, then $A$ is AP if and only if $-a_{33}a_{21} < -a_{23}a_{31}$. In particular, if we take $-a_{11}=a_{12}=a_{21}=-a_{23}=-a_{31}=a_{33}=1$, then $A$ is not AP. On the other hand, if we take $-a_{11}=a_{12}=a_{21}=-a_{23}=-a_{31}=1$ and $a_{33}=2$, then $A$ is AP.
	\end{itemize}
	
	\begin{minipage}{0.3\textwidth}
		\centering
		\includegraphics[scale=0.6]{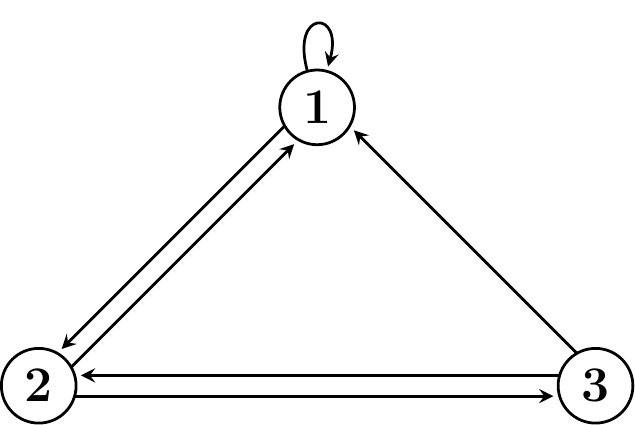}
		\captionof{Digraph}{}
		\label{D16}
	\end{minipage}
	\begin{minipage}{0.23\textwidth}
		\centering
		$\begin{bmatrix}
		* & + & 0 \\
		+ & 0 & + \\
		+ & + & 0
		\end{bmatrix}$
		\captionof{SP}{\textbf{RAP}}
	\end{minipage}
	\begin{minipage}{0.23\textwidth}
		\centering
		$\begin{bmatrix}
		* & + & 0 \\
		- & 0 & + \\
		+ & + & 0
		\end{bmatrix}$
		\captionof{SP}{\textbf{AAP}}
		\label{SP16.2}
	\end{minipage}
	\begin{minipage}{0.23\textwidth}
		\centering
		$\begin{bmatrix}
		* & + & 0 \\
		+ & 0 & + \\
		- & + & 0
		\end{bmatrix}$
		\captionof{SP}{\textbf{AAP}}
		\label{SP16.3}
	\end{minipage}
	\smallskip\\
	\begin{minipage}{0.3\textwidth}
		\centering
		\ 
	\end{minipage}
	\begin{minipage}{0.23\textwidth}
		\centering
		$\begin{bmatrix}
		* & + & 0 \\
		+ & 0 & + \\
		+ & - & 0
		\end{bmatrix}$
		\captionof{SP}{\textbf{AAP}}
		\label{SP16.4}
	\end{minipage}
	\begin{minipage}{0.23\textwidth}
		\centering
		$\begin{bmatrix}
		* & + & 0 \\
		- & 0 & + \\
		+ & - & 0
		\end{bmatrix}$				
		\captionof{SP}{\textbf{AAP}}
		\label{SP16.5}
	\end{minipage}
	\begin{minipage}{0.23\textwidth}
		\centering
		$\begin{bmatrix}
		+ & - & 0 \\
		- & 0 & + \\
		+ & + & 0
		\end{bmatrix}$			
		\captionof{SP}{\textbf{AAP}}
		\label{SP16.6}
	\end{minipage}
	\smallskip\\
	\begin{minipage}{0.3\textwidth}
		\centering
		\ 
	\end{minipage}
	\begin{minipage}{0.23\textwidth}
		\centering
		$\begin{bmatrix}
		+ & - & 0 \\
		+ & 0 & + \\
		- & + & 0
		\end{bmatrix}$
		\captionof{SP}{\textbf{AAP}}
		\label{SP16.7}
	\end{minipage}
	\begin{minipage}{0.23\textwidth}
		\centering
		$\begin{bmatrix}
		+ & - & 0 \\
		- & 0 & + \\
		- & + & 0
		\end{bmatrix}$				
		\captionof{SP}{\textbf{AAP}}
		\label{SP16.8}
	\end{minipage}
	\begin{minipage}{0.23\textwidth}
		\centering
		$\begin{bmatrix}
		- & - & 0 \\
		* & 0 & + \\
		* & * & 0
		\end{bmatrix}$			
		\captionof{SP}{\textbf{DNA}}
	\end{minipage}
	\smallskip\\
	\begin{minipage}{0.3\textwidth}
		\centering
		\ 
	\end{minipage}
	\begin{minipage}{0.23\textwidth}
		\centering
		$\begin{bmatrix}
		* & + & 0 \\
		+ & 0 & + \\
		- & - & 0
		\end{bmatrix}$
		\captionof{SP}{\textbf{DNA}}
	\end{minipage}
	\begin{minipage}{0.23\textwidth}
		\centering
		$\begin{bmatrix}
		* & + & 0 \\
		- & 0 & + \\
		- & + & 0
		\end{bmatrix}$				
		\captionof{SP}{\textbf{DNA}}
	\end{minipage}
	\begin{minipage}{0.23\textwidth}
		\centering
		$\begin{bmatrix}
		* & - & 0 \\
		+ & 0 & - \\
		+ & + & 0
		\end{bmatrix}$			
		\captionof{SP}{\textbf{DNA}}
	\end{minipage}
	\smallskip\\
	\begin{minipage}{0.53\textwidth}
		\centering
		\ 
	\end{minipage}
	\begin{minipage}{0.23\textwidth}
		\centering
		$\begin{bmatrix}
		+ & - & 0 \\
		+ & 0 & + \\
		+ & + & 0
		\end{bmatrix}$
		\captionof{SP}{\textbf{DNA}}
	\end{minipage}
	\begin{minipage}{0.23\textwidth}
		\centering
		$\begin{bmatrix}
		+ & - & 0 \\
		* & 0 & + \\
		* & - & 0
		\end{bmatrix}$			
		\captionof{SP}{\textbf{DNA}}
	\end{minipage}
	Suppose that the digraph of $A=[a_{ij}]\in M_3(\mathbb{R})$ is Digraph \ref{D16}. If the sign pattern of $A$ is
	\begin{itemize}
		\item SP \ref{SP16.2}, then $A$ is AP if and only if $-\dfrac{a_{12}a_{31}}{a_{32}}, -a_{11}-\dfrac{a_{21}a_{32}}{a_{31}} < -a_{11}-\dfrac{a_{23}a_{31}}{a_{21}}$. In particular, if we take $\pm a_{11}=a_{12}=-a_{21}=a_{23}=a_{31}=a_{32}=1$, then $A$ is not AP. On the other hand, if we take $\pm a_{11}=a_{12}=a_{23}=a_{31}=a_{32}=2$ and $a_{21}=-1$, then $A$ is AP.
		\item SP \ref{SP16.3}, then $A$ is AP if and only if $-\dfrac{a_{12}a_{31}}{a_{32}}, -a_{11}-\dfrac{a_{23}a_{31}}{a_{21}} < -a_{11}-\dfrac{a_{21}a_{32}}{a_{31}}$. In particular, if we take $\pm a_{11}=a_{12}=a_{21}=a_{23}=-a_{31}=a_{32}=1$, then $A$ is not AP. On the other hand, if we take $\pm a_{11}=a_{12}=a_{21}=a_{23}=a_{32}=2$ and $a_{21}=-1$, then $A$ is AP.
		\item SP \ref{SP16.4}, then $A$ is AP if and only if $-a_{11}-\dfrac{a_{21}a_{32}}{a_{31}} < -\dfrac{a_{12}a_{31}}{a_{32}}$. In particular, if we take $\pm a_{11}=a_{12}=a_{21}=a_{23}=-a_{31}=1$ and $a_{32}=-3$, then $A$ is not AP. On the other hand, if we take $\pm a_{11}=a_{12}=a_{21}=a_{23}=a_{31}=2$ and $a_{32}=-1$, then $A$ is AP.
		\item SP \ref{SP16.5}, then $A$ is AP if and only if $-a_{11}, -a_{11}-\dfrac{a_{21}a_{32}}{a_{31}} < -\dfrac{a_{12}a_{31}}{a_{32}}, -a_{11}-\dfrac{a_{23}a_{31}}{a_{21}}$. In particular, if we take $\pm a_{11}=a_{12}=-a_{21}=a_{23}=a_{31}=-a_{32}=1$, then $A$ is not AP. On the other hand, if we take $\pm a_{11}=a_{12}=a_{23}=a_{31}=2$ and $a_{21}=a_{32}=-1$, then $A$ is AP.
		\item SP \ref{SP16.6}, then $A$ is AP if and only if $-a_{11}, -a_{11}-\dfrac{a_{23}a_{31}}{a_{21}} < -\dfrac{a_{12}a_{31}}{a_{32}}, -a_{11}-\dfrac{a_{21}a_{32}}{a_{31}}$. In particular, if we take $a_{11}=-a_{12}=-a_{21}=a_{23}=a_{31}=a_{32}=1$, then $A$ is not AP. On the other hand, if we take $a_{11}=-a_{12}=-a_{21}=a_{31}=a_{32}=1$ and $a_{23}=0.1$, then $A$ is AP.
		\item SP \ref{SP16.7}, then $A$ is AP if and only if $-a_{11}, -a_{11}-\dfrac{a_{21}a_{32}}{a_{31}} < -\dfrac{a_{12}a_{31}}{a_{32}}, -a_{11}-\dfrac{a_{23}a_{31}}{a_{21}}$. In particular, if we take $a_{11}=-a_{12}=-a_{21}=a_{23}=a_{31}=a_{32}=1$, then $A$ is not AP. On the other hand, if we take $-a_{12}=a_{21}=a_{23}=-a_{31}=1$, $a_{11}=11$, and $a_{32}=0.1$, then $A$ is AP.
		\item SP \ref{SP16.8}, then $A$ is AP if and only if $-a_{11}, -a_{11}-\dfrac{a_{21}a_{32}}{a_{31}}, -a_{11}-\dfrac{a_{23}a_{31}}{a_{21}} < -\dfrac{a_{12}a_{31}}{a_{32}}$. In particular, if we take $a_{11}=-a_{12}=-a_{21}=a_{23}=-a_{31}=a_{32}=1$, then $A$ is not AP. On the other hand, if we take $-a_{12}=-a_{21}=a_{23}=-a_{31}=1$, $a_{11}=11$, and $a_{32}=0.1$, then $A$ is AP.
	\end{itemize}
	
	\begin{minipage}{0.3\textwidth}
		\centering
		\includegraphics[scale=0.6]{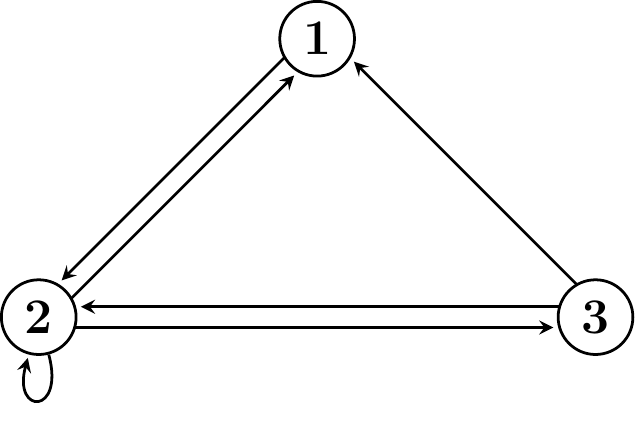}
		\captionof{Digraph}{}
		\label{D17}
	\end{minipage}
	\begin{minipage}{0.23\textwidth}
		\centering
		$\begin{bmatrix}
		0 & + & 0 \\
		+ & * & + \\
		+ & + & 0
		\end{bmatrix}$
		\captionof{SP}{\textbf{RAP}}
	\end{minipage}
	\begin{minipage}{0.23\textwidth}
		\centering
		$\begin{bmatrix}
		0 & + & 0 \\
		- & * & + \\
		+ & + & 0
		\end{bmatrix}$
		\captionof{SP}{\textbf{AAP}}
		\label{SP17.2}
	\end{minipage}
	\begin{minipage}{0.23\textwidth}
		\centering
		$\begin{bmatrix}
		0 & + & 0 \\
		* & * & - \\
		* & * & 0
		\end{bmatrix}$
		\captionof{SP}{\textbf{DNA}}
	\end{minipage}
	\smallskip\\
	\begin{minipage}{0.3\textwidth}
		\centering
		\ 
	\end{minipage}
	\begin{minipage}{0.23\textwidth}
		\centering
		$\begin{bmatrix}
		0 & + & 0 \\
		- & - & + \\
		+ & - & 0
		\end{bmatrix}$
		\captionof{SP}{\textbf{RAP}}
		\label{SP17.4}
	\end{minipage}
	\begin{minipage}{0.23\textwidth}
		\centering
		$\begin{bmatrix}
		0 & + & 0 \\
		+ & * & + \\
		- & + & 0
		\end{bmatrix}$				
		\captionof{SP}{\textbf{AAP}}
		\label{SP17.5}
	\end{minipage}
	\begin{minipage}{0.23\textwidth}
		\centering
		$\begin{bmatrix}
		0 & + & 0 \\
		- & * & + \\
		- & + & 0
		\end{bmatrix}$			
		\captionof{SP}{\textbf{DNA}}
	\end{minipage}
	\smallskip\\
	\begin{minipage}{0.53\textwidth}
		\centering
		\ 
	\end{minipage}
	\begin{minipage}{0.23\textwidth}
		\centering
		$\begin{bmatrix}
		0 & + & 0 \\
		- & + & + \\
		+ & - & 0
		\end{bmatrix}$				
		\captionof{SP}{\textbf{AAP}}
		\label{SP17.7}
	\end{minipage}
	\begin{minipage}{0.23\textwidth}
		\centering
		$\begin{bmatrix}
		0 & - & 0 \\
		+ & * & - \\
		+ & + & 0
		\end{bmatrix}$			
		\captionof{SP}{\textbf{DNA}}
	\end{minipage}
	
	Suppose that the digraph of $A=[a_{ij}]\in M_3(\mathbb{R})$ is Digraph \ref{D17}. If the sign pattern of $A$ is
	\begin{itemize}
		\item SP \ref{SP17.2}, then $A$ is AP if and only if $-a_{22}, -\dfrac{a_{21}a_{32}}{a_{31}}, -a_{33}-\dfrac{a_{12}a_{31}}{a_{32}} < -a_{33}-\dfrac{a_{23}a_{31}}{a_{21}}$. In particular, if we take $a_{12}=\pm a_{22}=a_{23}=a_{31}=a_{32}=1$ and $a_{21}=-3$, then $A$ is not AP. On the other hand, if we take $a_{12}=\pm a_{22}=a_{23}=a_{31}=a_{32}=2$ and $a_{21}=-1$, then $A$ is AP.
		\item SP \ref{SP17.4}, then take $p(x)=k_2x^2+k_1x+k_0$ such that $k_1,k_2>0$, $-a_{22}< \frac{k_1}{k_2}<-a_{22}-\max\{\frac{a_{23}a_{31}}{a_{21}},\frac{a_{31}a_{12}}{a_{32}}\}$ and $k_0$ be larger than all the diagonal entries of $-k_1A-k_2A^2$ so that $p(A)>0$.
		\item SP \ref{SP17.5}, then $A$ is AP if and only if $-a_{22}, -a_{33}-\dfrac{a_{23}a_{31}}{a_{21}}, -a_{33}-\dfrac{a_{12}a_{31}}{a_{32}} < -\dfrac{a_{21}a_{32}}{a_{31}}$. In particular, if we take $a_{12}=a_{21}=\pm a_{22}=a_{23}=a_{32}=1$ and $a_{31}=-3$, then $A$ is not AP. On the other hand, if we take $a_{12}=a_{21}=\pm a_{22}=a_{23}=a_{32}=2$ and $a_{31}=-1$, then $A$ is AP.
		\item SP \ref{SP17.7}, then $A$ is AP if and only if $-a_{22}, -\dfrac{a_{21}a_{32}}{a_{31}} < -a_{33}-\dfrac{a_{12}a_{31}}{a_{32}}, -a_{33}-\dfrac{a_{23}a_{31}}{a_{21}}$. In particular, if we take $a_{12}=-a_{21}=a_{23}=a_{31}=-a_{32}=1$ and $a_{22}=-5$, then $A$ is not AP. On the other hand, if we take $a_{12}=a_{22}=a_{23}=a_{31}=2$ and $a_{21}=a_{32}=-1$, then $A$ is AP.
	\end{itemize}
	\begin{minipage}{0.3\textwidth}
		\centering
		\includegraphics[scale=0.6]{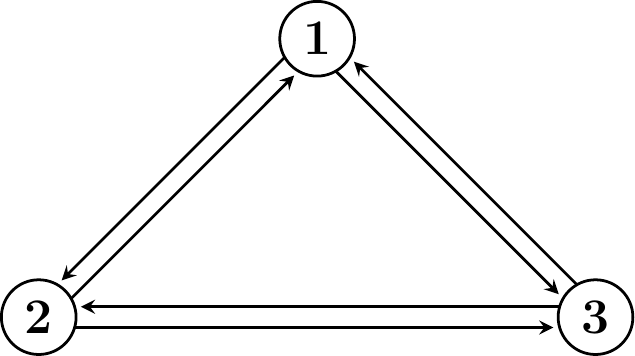}
		\captionof{Digraph}{}
		\label{D18}
	\end{minipage}
	\begin{minipage}{0.23\textwidth}
		\centering
		$\begin{bmatrix}
		0 & + & + \\
		+ & 0 & + \\
		+ & + & 0
		\end{bmatrix}$
		\captionof{SP}{\textbf{RAP}}
	\end{minipage}
	\begin{minipage}{0.23\textwidth}
		\centering
		$\begin{bmatrix}
		0 & - & + \\
		+ & 0 & + \\
		+ & + & 0
		\end{bmatrix}$
		\captionof{SP}{\textbf{AAP}}
		\label{SP18.2}
	\end{minipage}
	\begin{minipage}{0.23\textwidth}
		\centering
		$\begin{bmatrix}
		0 & - & + \\
		+ & 0 & + \\
		+ & - & 0
		\end{bmatrix}$
		\captionof{SP}{\textbf{DNA}}
	\end{minipage}
	\smallskip\\
	\begin{minipage}{0.3\textwidth}
		\centering
		\ 
	\end{minipage}
	\begin{minipage}{0.23\textwidth}
		\centering
		$\begin{bmatrix}
		0 & - & + \\
		+ & 0 & - \\
		- & + & 0
		\end{bmatrix}$
		\captionof{SP}{\textbf{RAP}}
		\label{SP18.4}
	\end{minipage}
	\begin{minipage}{0.23\textwidth}
		\centering
		$\begin{bmatrix}
		0 & - & + \\
		- & 0 & + \\
		+ & + & 0
		\end{bmatrix}$				
		\captionof{SP}{\textbf{AAP}}
		\label{SP18.5}
	\end{minipage}
	\begin{minipage}{0.23\textwidth}
		\centering
		$\begin{bmatrix}
		0 & - & - \\
		- & 0 & + \\
		+ & + & 0
		\end{bmatrix}$			
		\captionof{SP}{\textbf{DNA}}
	\end{minipage}
	\smallskip\\
	\begin{minipage}{0.53\textwidth}
		\centering
		\ 
	\end{minipage}
	\begin{minipage}{0.23\textwidth}
		\centering
		$\begin{bmatrix}
		0 & - & + \\
		+ & 0 & - \\
		+ & + & 0
		\end{bmatrix}$				
		\captionof{SP}{\textbf{AAP}}
		\label{SP18.7}
	\end{minipage}
	\begin{minipage}{0.23\textwidth}
		\centering
		$\begin{bmatrix}
		0 & - & - \\
		+ & 0 & - \\
		+ & + & 0
		\end{bmatrix}$			
		\captionof{SP}{\textbf{DNA}}
	\end{minipage}
	
	Suppose that the digraph of $A=[a_{ij}]\in M_3(\mathbb{R})$ is Digraph \ref{D18}. If the sign pattern of $A$ is
	\begin{itemize}
		\item SP \ref{SP18.2}, then $A$ is AP if and only if $-\dfrac{a_{12}a_{23}}{a_{13}}, -\dfrac{a_{23}a_{31}}{a_{21}}, -\dfrac{a_{21}a_{13}}{a_{23}}, -\dfrac{a_{21}a_{32}}{a_{31}}, -\dfrac{a_{12}a_{31}}{a_{32}} < -\dfrac{a_{13}a_{32}}{a_{12}}$. In particular, if we take $-a_{12}=a_{13}=a_{21}=a_{23}=a_{31}=a_{32}=1$, then $A$ is not AP. On the other hand, if we take $a_{13}=a_{21}=a_{23}=a_{31}=a_{32}=2$ and $a_{12}=-1$, then $A$ is AP.
		\item SP \ref{SP18.4}, then take $p(x)=x^2+k_0$ such that $k_0>\max\{-a_{12}a_{21}-a_{13}a_{31}-a_{12}a_{21}-a_{23}a_{32},-a_{32}a_{23}-a_{13}a_{31}\}$ so that $p(A)>0$.
		\item SP \ref{SP18.5}, then $A$ is AP if and only if $-\dfrac{a_{12}a_{23}}{a_{13}}, -\dfrac{a_{21}a_{13}}{a_{23}}, -\dfrac{a_{21}a_{32}}{a_{31}}, -\dfrac{a_{12}a_{31}}{a_{32}} < -\dfrac{a_{13}a_{32}}{a_{12}}, -\dfrac{a_{23}a_{31}}{a_{21}}$. In particular, if we take $-a_{12}=a_{13}=-a_{21}=a_{23}=a_{31}=a_{32}=1$, then $A$ is not AP. On the other hand, if we take $a_{13}=a_{23}=a_{31}=a_{32}=2$ and $a_{12}=a_{21}=-1$, then $A$ is AP.
		\item SP \ref{SP18.7}, then $A$ is AP if and only if $-\dfrac{a_{12}a_{23}}{a_{13}}, -\dfrac{a_{23}a_{31}}{a_{21}}, -\dfrac{a_{21}a_{32}}{a_{31}}, -\dfrac{a_{12}a_{31}}{a_{32}} < -\dfrac{a_{13}a_{32}}{a_{12}}, -\dfrac{a_{21}a_{13}}{a_{23}}$. In particular, if we take $-a_{12}=a_{13}=a_{21}=-a_{23}=a_{31}=a_{32}=1$, then $A$ is not AP. On the other hand, if we take $a_{13}=a_{21}=a_{31}=a_{32}=2$ and $a_{12}=a_{23}=-1$, then $A$ is AP.
	\end{itemize}
	
	\item \textbf{Graphs with 7 directed edges}. There are five nonequivalent strongly connected digraphs with 7 directed edges. We list down the nonequivalent sign pattern matrices corresponding to each digraph and except for SP \ref{SP19.4}, \ref{SP20.4}, and the AAP sign patterns, we use Theorems \ref{thm2}.2,\ref{thm2}.3 and \ref{thm4} to classify these sign patterns as indicated below. 
	
	\begin{minipage}{0.3\textwidth}
		\centering
		\includegraphics[scale=0.6]{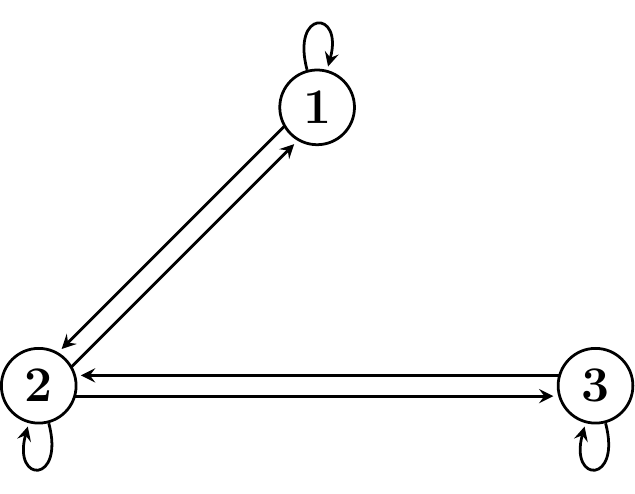}
		\captionof{Digraph}{}
		\label{D19}
	\end{minipage}
	\begin{minipage}{0.23\textwidth}
		\centering
		$\begin{bmatrix}
		* & + & 0 \\
		+ & * & + \\
		0 & + & *
		\end{bmatrix}$
		\captionof{SP}{\textbf{RAP}}
	\end{minipage}
	\begin{minipage}{0.23\textwidth}
		\centering
		$\begin{bmatrix}
		+ & - & 0 \\
		- & * & + \\
		0 & + & +
		\end{bmatrix}$
		\captionof{SP}{\textbf{AAP}}
		\label{SP19.2}
	\end{minipage}
	\begin{minipage}{0.23\textwidth}
		\centering
		$\begin{bmatrix}
		* & - & 0 \\
		+ & * & * \\
		0 & * & *
		\end{bmatrix}$
		\captionof{SP}{\textbf{DNA}}
	\end{minipage}
	\smallskip\\
	\begin{minipage}{0.3\textwidth}
		\centering
		\ 
	\end{minipage}
	\begin{minipage}{0.23\textwidth}
		\centering
		\ 
	\end{minipage}
	\begin{minipage}{0.23\textwidth}
		\centering
		$\begin{bmatrix}
		- & + & 0 \\
		+ & * & - \\
		0 & - & +
		\end{bmatrix}$
		\captionof{SP}{\textbf{RAP}}
		\label{SP19.4}
	\end{minipage}
	\begin{minipage}{0.23\textwidth}
		\centering
		$\begin{bmatrix}
		+ & + & 0 \\
		+ & * & - \\
		0 & - & -
		\end{bmatrix}$			
		\captionof{SP}{\textbf{DNA}}
	\end{minipage}
	
Suppose that the digraph of $A=[a_{ij}]\in M_3(\mathbb{R})$ is Digraph \ref{D19}. If the sign pattern of $A$ is
\begin{itemize}
\item SP \ref{SP19.2}, then $A$ is AP if and only if $a_{11}>a_{33}$. In particular, if we take $a_{11}=-a_{12}=-a_{21}=\pm a_{22}=a_{23}=a_{32}=a_{33}=1$, then $A$ is not AP. On the other hand, if we take $-a_{12}=-a_{21}=\pm a_{22}=a_{23}=a_{32}=a_{33}=1$ and $a_{11}=2$, then $A$ is AP.
\item SP \ref{SP19.4}, take $p(x)=k_2x^2+k_1x+k_0$ such that $k_2<0$, $-(a_{22}+a_{33})<\frac{k_1}{k_2}< -(a_{11}+a_{22})$ and $k_0$ larger than all the diagonal entries of $-k_1A-k_2A^2$ so that $p(A)>0$. 

\end{itemize}
	\begin{minipage}{0.3\textwidth}
		\centering
		\includegraphics[scale=0.6]{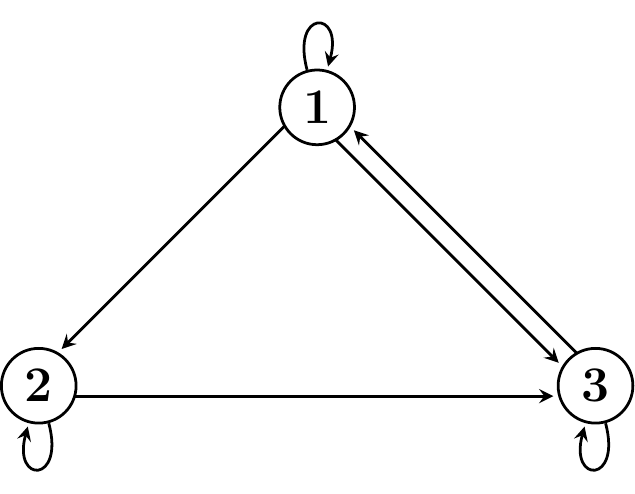}
		\captionof{Digraph}{}
		\label{D20}
	\end{minipage}
	\begin{minipage}{0.23\textwidth}
		\centering
		$\begin{bmatrix}
		* & + & + \\
		0 & * & + \\
		+ & 0 & *
		\end{bmatrix}$
		\captionof{SP}{\textbf{RAP}}
	\end{minipage}
	\begin{minipage}{0.23\textwidth}
		\centering
		$\begin{bmatrix}
		+ & + & - \\
		0 & + & + \\
		+ & 0 & *
		\end{bmatrix}$
		\captionof{SP}{\textbf{AAP}}
		\label{SP20.2}
	\end{minipage}
	\begin{minipage}{0.23\textwidth}
		\centering
		$\begin{bmatrix}
		* & + & * \\
		0 & * & - \\
		* & 0 & *
		\end{bmatrix}$
		\captionof{SP}{\textbf{DNA}}
	\end{minipage}
	\smallskip\\
	\begin{minipage}{0.3\textwidth}
		\centering
		\ 
	\end{minipage}
	\begin{minipage}{0.23\textwidth}
		\centering
		$\begin{bmatrix}
		- & + & - \\
		0 & + & + \\
		+ & 0 & -
		\end{bmatrix}$
		\captionof{SP}{\textbf{RAP}}
		\label{SP20.4}
	\end{minipage}
	\begin{minipage}{0.23\textwidth}
		\centering
		$\begin{bmatrix}
		* & + & - \\
		0 & - & + \\
		+ & 0 & *
		\end{bmatrix}$				
		\captionof{SP}{\textbf{AAP}}
		\label{SP20.5}
	\end{minipage}
	\begin{minipage}{0.23\textwidth}
		\centering
		$\begin{bmatrix}
		* & + & + \\
		0 & * & + \\
		- & 0 & *
		\end{bmatrix}$			
		\captionof{SP}{\textbf{DNA}}
	\end{minipage}
	\smallskip\\
	\begin{minipage}{0.3\textwidth}
		\centering
		\ 
	\end{minipage}
	\begin{minipage}{0.23\textwidth}
		\centering
		$\begin{bmatrix}
		+ & + & - \\
		0 & + & + \\
		- & 0 & +
		\end{bmatrix}$
		\captionof{SP}{\textbf{AAP}}
		\label{SP20.7}
	\end{minipage}
	\begin{minipage}{0.23\textwidth}
		\centering
		$\begin{bmatrix}
		* & - & + \\
		0 & + & - \\
		+ & 0 & *
		\end{bmatrix}$				
		\captionof{SP}{\textbf{AAP}}
		\label{SP20.8}
	\end{minipage}
	\begin{minipage}{0.23\textwidth}
		\centering
		$\begin{bmatrix}
		- & + & - \\
		0 & + & + \\
		- & 0 & *
		\end{bmatrix}$			
		\captionof{SP}{\textbf{DNA}}
	\end{minipage}

	Suppose that the digraph of $A=[a_{ij}]\in M_3(\mathbb{R})$ is Digraph \ref{D20}. If the sign pattern of $A$ is
	\begin{itemize}
		\item SP \ref{SP20.2}, then $A$ is AP if and only if $-a_{11}-a_{22}, -a_{22}-a_{33}, -a_{11}-a_{33} < -a_{11}-a_{33}-\dfrac{a_{12}a_{23}}{a_{13}}$. In particular, if we take $a_{12}=a_{22}=a_{23}=a_{31}=1$, $a_{11}=\pm a_{33}=2$, and $a_{13}=-5$, then $A$ is not AP. On the other hand, if we take $a_{11}=a_{12}=-a_{13}=a_{23}=a_{31}=\pm a_{33}=1$ and $a_{22}=3$, then $A$ is AP.
		\item SP \ref{SP20.4}, then take $p(x)=k_2x^2+k_1x+k_0$ such that $k_1,k_2>0$, $-(a_{11}+a_{33})\leq \frac{k_1}{k_2}\leq -(a_{11}+a_{33})-\frac{a_{12}a_{23}}{a_{13}}$ and $k_0$ is larger than all the diagonal entries of $-k_1A-k_2A^2$ so that $p(A)>0$.
		\item SP \ref{SP20.5}, then $A$ is AP if and only if $-a_{11}-a_{22}, -a_{22}-a_{33}, -a_{11}-a_{33} < -a_{11}-a_{33}-\dfrac{a_{12}a_{23}}{a_{13}}$. In particular, if we take $\pm a_{11}=a_{12}=-a_{13}=a_{23}=a_{31}=\pm a_{33}=1$ and $a_{22}=-10$, then $A$ is not AP. On the other hand, if we take $\pm a_{11}=a_{12}=-a_{22}=a_{23}=a_{31}=\pm a_{33}=1$ and $a_{13}=-0.1$, then $A$ is AP.
		\item SP \ref{SP20.7}, then $A$ is AP if and only if $-a_{11}-a_{33}, -a_{11}-a_{33}-\dfrac{a_{12}a_{23}}{a_{13}} < -a_{11}-a_{22}, -a_{22}-a_{33}$. In particular, if we take $a_{11}=a_{12}=-a_{13}=a_{22}=a_{23}=-a_{31}=a_{33}=1$, then $A$ is not AP. On the other hand, if we take $\pm a_{12}=a_{22}=a_{23}=-a_{31}=1$ and $a_{11}=-a_{13}=a_{33}=5$, then $A$ is AP.
		\item SP \ref{SP20.8}, then $A$ is AP if and only if $-a_{11}-a_{22}, -a_{22}-a_{33} < -a_{11}-a_{33}, -a_{11}-a_{33}-\dfrac{a_{12}a_{23}}{a_{13}}$. In particular, if we take $\pm a_{11}=-a_{12}=a_{22}=-a_{23}=a_{31}=\pm a_{33}=1$ and $a_{13}=0.1$, then $A$ is not AP. On the other hand, if we take $\pm a_{11}=-a_{12}=a_{13}=-a_{23}=a_{31}=\pm a_{33}=1$ and $a_{22}=5$, then $A$ is AP.
	\end{itemize}
	
	\begin{minipage}{0.3\textwidth}
		\centering
		\includegraphics[scale=0.6]{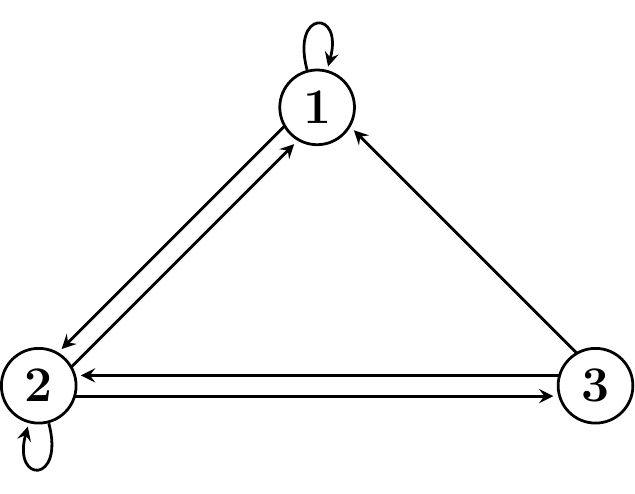}
		\captionof{Digraph}{}
		\label{D21}
	\end{minipage}
	\begin{minipage}{0.23\textwidth}
		\centering
		$\begin{bmatrix}
		* & + & 0 \\
		+ & * & + \\
		+ & + & 0
		\end{bmatrix}$
		\captionof{SP}{\textbf{RAP}}
	\end{minipage}
	\begin{minipage}{0.23\textwidth}
		\centering
		$\begin{bmatrix}
		* & + & 0 \\
		- & * & + \\
		+ & + & 0
		\end{bmatrix}$
		\captionof{SP}{\textbf{AAP}}
		\label{SP21.2}
	\end{minipage}
	\begin{minipage}{0.23\textwidth}
		\centering
		$\begin{bmatrix}
		* & + & 0 \\
		+ & * & + \\
		+ & - & 0
		\end{bmatrix}$
		\captionof{SP}{\textbf{AAP}}
		\label{SP21.3}
	\end{minipage}
	\smallskip\\
	\begin{minipage}{0.3\textwidth}
		\centering
		\ 
	\end{minipage}
	\begin{minipage}{0.23\textwidth}
		\centering
		$\begin{bmatrix}
		* & + & 0 \\
		+ & * & + \\
		- & + & 0
		\end{bmatrix}$
		\captionof{SP}{\textbf{AAP}}
		\label{SP21.4}
	\end{minipage}
	\begin{minipage}{0.23\textwidth}
		\centering
		$\begin{bmatrix}
		+ & - & 0 \\
		- & * & + \\
		+ & + & 0
		\end{bmatrix}$				
		\captionof{SP}{\textbf{AAP}}
		\label{SP21.5}
	\end{minipage}
	\begin{minipage}{0.23\textwidth}
		\centering
		$\begin{bmatrix}
		+ & - & 0 \\
		+ & * & + \\
		- & + & 0
		\end{bmatrix}$			
		\captionof{SP}{\textbf{AAP}}
		\label{SP21.6}
	\end{minipage}
	\smallskip\\
	\begin{minipage}{0.3\textwidth}
		\centering
		\ 
	\end{minipage}
	\begin{minipage}{0.23\textwidth}
		\centering
		$\begin{bmatrix}
		* & + & 0 \\
		- & * & + \\
		+ & - & 0
		\end{bmatrix}$
		\captionof{SP}{\textbf{AAP}}
		\label{SP21.7}
	\end{minipage}
	\begin{minipage}{0.23\textwidth}
		\centering
		$\begin{bmatrix}
		- & + & 0 \\
		+ & - & - \\
		- & + & 0
		\end{bmatrix}$				
		\captionof{SP}{\textbf{AAP}}
		\label{SP21.8}
	\end{minipage}
	\begin{minipage}{0.23\textwidth}
		\centering
		$\begin{bmatrix}
		- & + & 0 \\
		+ & * & - \\
		+ & - & 0
		\end{bmatrix}$			
		\captionof{SP}{\textbf{AAP}}
		\label{SP21.9}
	\end{minipage}
	\begin{minipage}{0.3\textwidth}
		\centering
		\ 
	\end{minipage}
	\begin{minipage}{0.23\textwidth}
		\centering
		$\begin{bmatrix}
		* & - & 0 \\
		+ & * & + \\
		+ & + & 0
		\end{bmatrix}$
		\captionof{SP}{\textbf{DNA}}
	\end{minipage}
	\begin{minipage}{0.23\textwidth}
		\centering
		$\begin{bmatrix}
		* & + & 0 \\
		+ & * & - \\
		+ & + & 0
		\end{bmatrix}$				
		\captionof{SP}{\textbf{DNA}}
	\end{minipage}
	\begin{minipage}{0.23\textwidth}
		\centering
		$\begin{bmatrix}
		* & - & 0 \\
		+ & * & - \\
		+ & + & 0
		\end{bmatrix}$			
		\captionof{SP}{\textbf{DNA}}
	\end{minipage}
	\begin{minipage}{0.3\textwidth}
		\centering
		\ 
	\end{minipage}
	\begin{minipage}{0.23\textwidth}
		\centering
		$\begin{bmatrix}
		* & - & 0 \\
		+ & * & + \\
		+ & - & 0
		\end{bmatrix}$
		\captionof{SP}{\textbf{DNA}}
	\end{minipage}
	\begin{minipage}{0.23\textwidth}
		\centering
		$\begin{bmatrix}
		* & + & 0 \\
		- & * & - \\
		+ & + & 0
		\end{bmatrix}$				
		\captionof{SP}{\textbf{DNA}}
	\end{minipage}
	\begin{minipage}{0.23\textwidth}
		\centering
		$\begin{bmatrix}
		* & + & 0 \\
		- & * & + \\
		- & + & 0
		\end{bmatrix}$			
		\captionof{SP}{\textbf{DNA}}
	\end{minipage}
	\begin{minipage}{0.3\textwidth}
		\centering
		\ 
	\end{minipage}
	\begin{minipage}{0.23\textwidth}
		\centering
		$\begin{bmatrix}
		* & + & 0 \\
		+ & * & + \\
		- & - & 0
		\end{bmatrix}$
		\captionof{SP}{\textbf{DNA}}
	\end{minipage}
	\begin{minipage}{0.23\textwidth}
		\centering
		$\begin{bmatrix}
		- & - & 0 \\
		- & * & + \\
		+ & + & 0
		\end{bmatrix}$				
		\captionof{SP}{\textbf{DNA}}
	\end{minipage}
	\begin{minipage}{0.23\textwidth}
		\centering
		$\begin{bmatrix}
		- & - & 0 \\
		+ & * & + \\
		- & + & 0
		\end{bmatrix}$			
		\captionof{SP}{\textbf{DNA}}
	\end{minipage}
	\begin{minipage}{0.3\textwidth}
		\centering
		\ 
	\end{minipage}
	\begin{minipage}{0.23\textwidth}
		\centering
		$\begin{bmatrix}
		+ & + & 0 \\
		+ & * & - \\
		- & + & 0
		\end{bmatrix}$
		\captionof{SP}{\textbf{DNA}}
	\end{minipage}
	\begin{minipage}{0.23\textwidth}
		\centering
		$\begin{bmatrix}
		- & + & 0 \\
		+ & + & - \\
		- & + & 0
		\end{bmatrix}$				
		\captionof{SP}{\textbf{DNA}}
	\end{minipage}
	\begin{minipage}{0.23\textwidth}
		\centering
		$\begin{bmatrix}
		+ & + & 0 \\
		+ & * & - \\
		+ & - & 0
		\end{bmatrix}$			
		\captionof{SP}{\textbf{DNA}}
	\end{minipage}
	
	Suppose that the digraph of $A=[a_{ij}]\in M_3(\mathbb{R})$ is Digraph \ref{D21}. If the sign pattern of $A$ is
	\begin{itemize}
		\item SP \ref{SP21.2}, then $A$ is AP if and only if $-a_{11}-a_{22}, -a_{22}, -a_{11}-\dfrac{a_{21}a_{32}}{a_{31}}, -a_{22}-\dfrac{a_{12}a_{31}}{a_{32}} < -a_{11}-a_{22}-\dfrac{a_{23}a_{31}}{a_{21}}$. In particular, if we take $\pm a_{11}=a_{12}=\pm a_{22}=a_{23}=a_{31}=a_{32}=1$ and $a_{21}=-10$, then $A$ is not AP. On the other hand, if we take $\pm a_{11}=a_{12}=\pm a_{22}=a_{23}=a_{31}=a_{32}=1$ and $a_{21}=-0.1$, then $A$ is AP.
		\item SP \ref{SP21.3}, then $A$ is AP if and only if $-a_{11}-a_{22}, -a_{22}, -a_{11}-\dfrac{a_{21}a_{32}}{a_{31}}, -a_{11}-a_{22}-\dfrac{a_{23}a_{31}}{a_{21}} < -a_{22}-\dfrac{a_{12}a_{31}}{a_{32}}$. In particular, if we take $\pm a_{11}=a_{12}=\pm a_{22}=a_{23}=a_{31}=-a_{32}=1$ and $a_{21}=-10$, then $A$ is not AP. On the other hand, if we take $\pm a_{11}=a_{12}=a_{21}=\pm a_{22}=a_{23}=a_{31}=1$ and $a_{32}=-0.1$, then $A$ is AP.
		\item SP \ref{SP21.4}, then $A$ is AP if and only if $-a_{11}-a_{22}, -a_{22}, -a_{22}-\dfrac{a_{12}a_{31}}{a_{32}}, -a_{11}-a_{22}-\dfrac{a_{23}a_{31}}{a_{21}} < -a_{11}-\dfrac{a_{21}a_{32}}{a_{31}}$. In particular, if we take $\pm a_{11}=a_{12}=a_{21}=\pm a_{22}=a_{23}=a_{32}=1$ and $a_{31}=-10$, then $A$ is not AP. On the other hand, if we take $\pm a_{11}=a_{12}=a_{21}=\pm a_{22}=a_{23}=a_{32}=1$ and $a_{31}=-0.1$, then $A$ is AP.
		\item SP \ref{SP21.5}, then $A$ is AP if and only if $-a_{11}-a_{22}, -a_{11}-a_{22}-\dfrac{a_{23}a_{31}}{a_{21}} < -a_{22}, -a_{22}-\dfrac{a_{12}a_{31}}{a_{32}}, -a_{11}-\dfrac{a_{21}a_{32}}{a_{31}}$. In particular, if we take $a_{11}=-a_{12}=\pm a_{22}=a_{23}=a_{31}=a_{32}=1$ and $a_{21}=-0.1$, then $A$ is not AP. On the other hand, if we take $a_{11}=-a_{12}=\pm a_{22}=a_{23}=a_{31}=a_{32}=1$ and $a_{21}=-10$, then $A$ is AP.
		\item SP \ref{SP21.6}, then $A$ is AP if and only if $-a_{11}-a_{22}, -a_{11}-\dfrac{a_{21}a_{32}}{a_{31}} < -a_{22}, -a_{22}-\dfrac{a_{12}a_{31}}{a_{32}}, -a_{11}-a_{22}-\dfrac{a_{23}a_{31}}{a_{21}}$. In particular, if we take $a_{11}=-a_{12}=a_{21}=\pm a_{22}=a_{23}=a_{32}=1$ and $a_{31}=-0.1$, then $A$ is not AP. On the other hand, if we take $-a_{12}=a_{21}=\pm a_{22}=a_{23}=a_{32}=1$, $a_{11}=10$, and $a_{31}=-5$, then $A$ is AP.
		\item SP \ref{SP21.7}, then $A$ is AP if and only if $-a_{22}, -a_{11}-a_{22}, -a_{11}-\dfrac{a_{21}a_{32}}{a_{31}} < -a_{22}-\dfrac{a_{12}a_{31}}{a_{32}}, -a_{11}-a_{22}-\dfrac{a_{23}a_{31}}{a_{21}}$. In particular, if we take $\pm a_{11}=a_{12}=-a_{21}=\pm a_{22}=a_{23}=-a_{32}=1$ and $a_{31}=0.1$, then $A$ is not AP. On the other hand, if we take $\pm a_{11}=a_{12}=-a_{21}=\pm a_{22}=a_{23}=-a_{32}=1$ and $a_{31}=10$, then $A$ is AP.
		\item SP \ref{SP21.8}, then $A$ is AP if and only if $-a_{22},  -a_{11}-\dfrac{a_{21}a_{32}}{a_{31}} < -a_{11}-a_{22}, -a_{22}-\dfrac{a_{12}a_{31}}{a_{32}}, -a_{11}-a_{22}-\dfrac{a_{23}a_{31}}{a_{21}}$. In particular, if we take $-a_{11}=a_{21}=-a_{22}=-a_{23}=-a_{31}=a_{32}=1$ and $a_{12}=10$, then $A$ is not AP. On the other hand, if we take $a_{21}=-a_{23}=-a_{31}=a_{32}=1$, $a_{12}=-a_{22}=10$, and $a_{11}=-3$, then $A$ is AP.
		\item SP \ref{SP21.8}, then $A$ is AP if and only if $-a_{22}, -a_{22}-\dfrac{a_{12}a_{31}}{a_{32}} < -a_{11}-a_{22}, -a_{11}-\dfrac{a_{21}a_{32}}{a_{31}}, -a_{11}-a_{22}-\dfrac{a_{23}a_{31}}{a_{21}}$. In particular, if we take $-a_{11}=a_{12}=a_{21}=\pm a_{22}=-a_{23}=a_{31}=1$ and $a_{32}=-0.1$, then $A$ is not AP. On the other hand, if we take $-a_{11}=a_{12}=a_{21}=\pm a_{22}=-a_{23}=a_{31}=1$ and $a_{32}=-10$, then $A$ is AP.
	\end{itemize}
	
	\begin{minipage}{0.3\textwidth}
		\centering
		\includegraphics[scale=0.6]{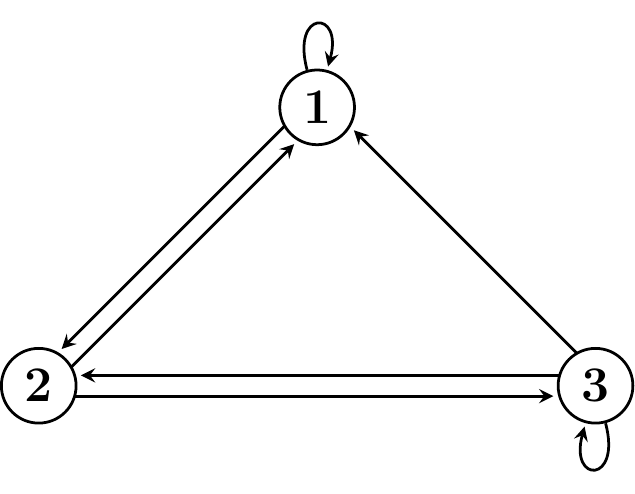}
		\captionof{Digraph}{}
		\label{D22}
	\end{minipage}
	\begin{minipage}{0.23\textwidth}
		\centering
		$\begin{bmatrix}
		* & + & 0 \\
		+ & 0 & + \\
		+ & + & *
		\end{bmatrix}$
		\captionof{SP}{\textbf{RAP}}
	\end{minipage}
	\begin{minipage}{0.23\textwidth}
		\centering
		$\begin{bmatrix}
		* & + & 0 \\
		- & 0 & + \\
		+ & + & *
		\end{bmatrix}$
		\captionof{SP}{\textbf{AAP}}
		\label{SP22.2}
	\end{minipage}
	\begin{minipage}{0.23\textwidth}
		\centering
		$\begin{bmatrix}
		* & + & 0 \\
		+ & 0 & + \\
		- & + & *
		\end{bmatrix}$
		\captionof{SP}{\textbf{AAP}}
		\label{SP22.3}
	\end{minipage}
	\smallskip\\
	\begin{minipage}{0.3\textwidth}
		\centering
		\ 
	\end{minipage}
	\begin{minipage}{0.23\textwidth}
		\centering
		$\begin{bmatrix}
		- & - & 0 \\
		- & 0 & + \\
		+ & + & -
		\end{bmatrix}$
		\captionof{SP}{\textbf{AAP}}
		\label{SP22.4}
	\end{minipage}
	\begin{minipage}{0.23\textwidth}
		\centering
		$\begin{bmatrix}
		+ & - & 0 \\
		- & 0 & + \\
		+ & + & *
		\end{bmatrix}$				
		\captionof{SP}{\textbf{AAP}}
		\label{SP22.5}
	\end{minipage}
	\begin{minipage}{0.23\textwidth}
		\centering
		$\begin{bmatrix}
		+ & - & 0 \\
		+ & 0 & + \\
		- & + & *
		\end{bmatrix}$			
		\captionof{SP}{\textbf{AAP}}
		\label{SP22.6}
	\end{minipage}
	\smallskip\\
	\begin{minipage}{0.3\textwidth}
		\centering
		\ 
	\end{minipage}
	\begin{minipage}{0.23\textwidth}
		\centering
		$\begin{bmatrix}
		* & + & 0 \\
		- & 0 & + \\
		+ & - & *
		\end{bmatrix}$
		\captionof{SP}{\textbf{AAP}}
		\label{SP22.7}
	\end{minipage}
	\begin{minipage}{0.23\textwidth}
		\centering
		$\begin{bmatrix}
		* & - & 0 \\
		+ & 0 & + \\
		+ & + & *
		\end{bmatrix}$				
		\captionof{SP}{\textbf{DNA}}
	\end{minipage}
	\begin{minipage}{0.23\textwidth}
		\centering
		$\begin{bmatrix}
		* & - & 0 \\
		+ & 0 & - \\
		+ & + & *
		\end{bmatrix}$			
		\captionof{SP}{\textbf{DNA}}
	\end{minipage}
	\begin{minipage}{0.3\textwidth}
		\centering
		\ 
	\end{minipage}
	\begin{minipage}{0.23\textwidth}
		\centering
		$\begin{bmatrix}
		* & - & 0 \\
		+ & 0 & + \\
		+ & - & *
		\end{bmatrix}$
		\captionof{SP}{\textbf{DNA}}
	\end{minipage}
	\begin{minipage}{0.23\textwidth}
		\centering
		$\begin{bmatrix}
		* & + & 0 \\
		- & 0 & + \\
		- & + & *
		\end{bmatrix}$				
		\captionof{SP}{\textbf{DNA}}
	\end{minipage}
	\begin{minipage}{0.23\textwidth}
		\centering
		$\begin{bmatrix}
		- & - & 0 \\
		- & 0 & + \\
		+ & + & +
		\end{bmatrix}$			
		\captionof{SP}{\textbf{DNA}}
	\end{minipage}
	\begin{minipage}{0.76\textwidth}
		\centering
		\ 
	\end{minipage}
	\begin{minipage}{0.23\textwidth}
		\centering
		$\begin{bmatrix}
		- & - & 0 \\
		+ & 0 & + \\
		- & + & *
		\end{bmatrix}$			
		\captionof{SP}{\textbf{DNA}}
	\end{minipage}
	
	Suppose that the digraph of $A=[a_{ij}]\in M_3(\mathbb{R})$ is Digraph \ref{D22}. If the sign pattern of $A$ is
	\begin{itemize}
		\item SP \ref{SP22.2}, then $A$ is AP if and only if $-a_{11}, -a_{33}, -a_{11}-a_{33}-\dfrac{a_{21}a_{32}}{a_{31}}, -a_{33}-\dfrac{a_{12}a_{31}}{a_{32}} < -a_{11}-\dfrac{a_{23}a_{31}}{a_{21}}$. In particular, if we take $\pm a_{11}=a_{12}=a_{23}=a_{31}=a_{32}=\pm a_{33}=1$ and $a_{21}=-10$, then $A$ is not AP. On the other hand, if we take $\pm a_{11}=a_{12}=a_{23}=a_{31}=a_{32}=\pm a_{33}=1$ and $a_{21}=-0.1$, then $A$ is AP.
		\item SP \ref{SP22.3}, then $A$ is AP if and only if $-a_{11}, -a_{33}, -a_{11}-\dfrac{a_{23}a_{31}}{a_{21}}, -a_{33}-\dfrac{a_{12}a_{31}}{a_{32}} < -a_{11}-a_{33}-\dfrac{a_{21}a_{32}}{a_{31}}$. In particular, if we take $\pm a_{11}=a_{12}=a_{21}=a_{23}=a_{32}=\pm a_{33}=1$ and $a_{31}=-10$, then $A$ is not AP. On the other hand, if we take $\pm a_{11}=a_{12}=a_{21}=a_{23}=a_{32}=\pm a_{33}=1$ and $a_{31}=-0.1$, then $A$ is AP.
		\item SP \ref{SP22.4}, then $A$ is AP if and only if $-a_{11},  -a_{11}-\dfrac{a_{23}a_{31}}{a_{21}} < -a_{33}, -a_{11}-a_{33}-\dfrac{a_{21}a_{32}}{a_{31}}, -a_{33}-\dfrac{a_{12}a_{31}}{a_{32}}$. In particular, if we take $-a_{11}=-a_{12}=-a_{21}=a_{23}=a_{31}=-a_{33}=1$ and $a_{32}=10$, then $A$ is not AP. On the other hand, if we take $-a_{11}=-a_{12}=-a_{21}=a_{23}=a_{31}=1$, $a_{32}=0.1$, and $a_{33}=-10$, then $A$ is AP.
		\item SP \ref{SP22.5}, then $A$ is AP if and only if $-a_{11},  -a_{11}-\dfrac{a_{23}a_{31}}{a_{21}} < -a_{33}, -a_{11}-a_{33}-\dfrac{a_{21}a_{32}}{a_{31}}, -a_{33}-\dfrac{a_{12}a_{31}}{a_{32}}$. In particular, if we take $a_{11}=-a_{12}=-a_{21}=a_{31}=\pm a_{33}=1$ and $a_{23}=a_{32}=10$, then $A$ is not AP. On the other hand, if we take $-a_{12}=-a_{21}=a_{23}=a_{31}=\pm a_{33}=1$, $a_{11}=15$, and $a_{32}=10$, then $A$ is AP.
		\item SP \ref{SP22.6}, then $A$ is AP if and only if $-a_{11}, -a_{11}-a_{33}-\dfrac{a_{21}a_{32}}{a_{31}} < -a_{33}, -a_{11}-\dfrac{a_{23}a_{31}}{a_{21}}, -a_{33}-\dfrac{a_{12}a_{31}}{a_{32}}$. In particular, if we take $a_{11}=-a_{12}=a_{21}=a_{23}=-a_{31}=\pm a_{33}=1$ and $a_{32}=10$, then $A$ is not AP. On the other hand, if we take $a_{21}=a_{23}=\pm a_{33}=1$, $a_{11}=15$, $a_{31}=-5$ and $-a_{12}=a_{32}=0.1$, then $A$ is AP.
		\item SP \ref{SP22.7}, then $A$ is AP if and only if $-a_{11}, -a_{33}, -a_{11}-a_{33}-\dfrac{a_{21}a_{32}}{a_{31}} < -a_{11}-\dfrac{a_{23}a_{31}}{a_{21}}, -a_{33}-\dfrac{a_{12}a_{31}}{a_{32}}$. In particular, if we take $\pm a_{11}=a_{12}=-a_{21}=a_{23}=a_{31}=1$ and $-a_{32}=\pm a_{33}=10$, then $A$ is not AP. On the other hand, if we take $\pm a_{11}=a_{12}=a_{23}=a_{31}=\pm a_{33}=1$ and $a_{23}=a_{32}=-0.1$, then $A$ is AP.
	\end{itemize}
	
	\begin{minipage}{0.3\textwidth}
		\centering
		\includegraphics[scale=0.6]{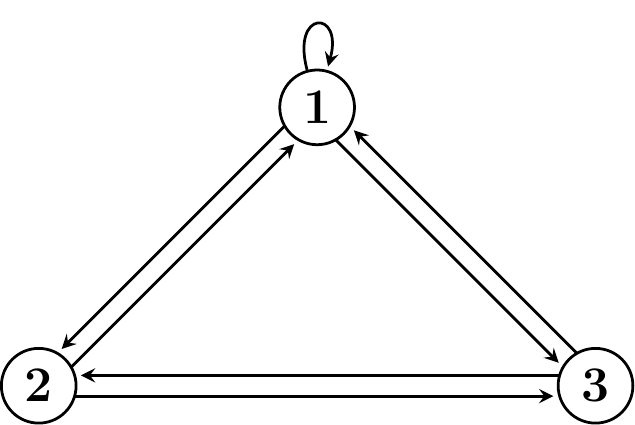}
		\captionof{Digraph}{}
		\label{D23}
	\end{minipage}
	\begin{minipage}{0.23\textwidth}
		\centering
		$\begin{bmatrix}
		* & + & + \\
		+ & 0 & + \\
		+ & + & 0
		\end{bmatrix}$
		\captionof{SP}{\textbf{RAP}}
	\end{minipage}
	\begin{minipage}{0.23\textwidth}
		\centering
		$\begin{bmatrix}
		* & - & + \\
		+ & 0 & + \\
		+ & + & 0
		\end{bmatrix}$
		\captionof{SP}{\textbf{AAP}}
		\label{SP23.2}
	\end{minipage}
	\begin{minipage}{0.23\textwidth}
		\centering
		$\begin{bmatrix}
		* & + & + \\
		+ & 0 & - \\
		+ & + & 0
		\end{bmatrix}$
		\captionof{SP}{\textbf{AAP}}
		\label{SP23.3}
	\end{minipage}
	\smallskip\\
	\begin{minipage}{0.3\textwidth}
		\centering
		\ 
	\end{minipage}
	\begin{minipage}{0.23\textwidth}
		\centering
		$\begin{bmatrix}
		* & - & + \\
		- & 0 & + \\
		+ & + & 0
		\end{bmatrix}$
		\captionof{SP}{\textbf{AAP}}
		\label{SP23.4}
	\end{minipage}
	\begin{minipage}{0.23\textwidth}
		\centering
		$\begin{bmatrix}
		* & + & + \\
		+ & 0 & - \\
		+ & - & 0
		\end{bmatrix}$				
		\captionof{SP}{\textbf{AAP}}
		\label{SP23.5}
	\end{minipage}
	\begin{minipage}{0.23\textwidth}
		\centering
		$\begin{bmatrix}
		* & - & + \\
		+ & 0 & - \\
		+ & + & 0
		\end{bmatrix}$			
		\captionof{SP}{\textbf{AAP}}
		\label{SP23.6}
	\end{minipage}
	\smallskip\\
	\begin{minipage}{0.3\textwidth}
		\centering
		\ 
	\end{minipage}
	\begin{minipage}{0.23\textwidth}
		\centering
		$\begin{bmatrix}
		* & - & + \\
		+ & 0 & + \\
		- & + & 0
		\end{bmatrix}$
		\captionof{SP}{\textbf{AAP}}
		\label{SP23.7}
	\end{minipage}
	\begin{minipage}{0.23\textwidth}
		\centering
		$\begin{bmatrix}
		+ & - & - \\
		- & 0 & + \\
		+ & + & 0
		\end{bmatrix}$				
		\captionof{SP}{\textbf{AAP}}
		\label{SP23.8}
	\end{minipage}
	\begin{minipage}{0.23\textwidth}
		\centering
		$\begin{bmatrix}
		* & - & + \\
		+ & 0 & - \\
		- & + & 0
		\end{bmatrix}$			
		\captionof{SP}{\textbf{AAP}}
		\label{SP23.9}
	\end{minipage}
	\begin{minipage}{0.3\textwidth}
		\centering
		\ 
	\end{minipage}
	\begin{minipage}{0.23\textwidth}
		\centering
		$\begin{bmatrix}
		* & - & - \\
		+ & 0 & + \\
		+ & + & 0
		\end{bmatrix}$
		\captionof{SP}{\textbf{DNA}}
	\end{minipage}
	\begin{minipage}{0.23\textwidth}
		\centering
		$\begin{bmatrix}
		* & - & + \\
		+ & 0 & + \\
		+ & - & 0
		\end{bmatrix}$				
		\captionof{SP}{\textbf{DNA}}
	\end{minipage}
	\begin{minipage}{0.23\textwidth}
		\centering
		$\begin{bmatrix}
		* & - & + \\
		- & 0 & - \\
		+ & + & 0
		\end{bmatrix}$			
		\captionof{SP}{\textbf{DNA}}
	\end{minipage}
	\begin{minipage}{0.3\textwidth}
		\centering
		\ 
	\end{minipage}
	\begin{minipage}{0.23\textwidth}
		\centering
		$\begin{bmatrix}
		* & - & - \\
		+ & 0 & - \\
		+ & + & 0
		\end{bmatrix}$
		\captionof{SP}{\textbf{DNA}}
	\end{minipage}
	\begin{minipage}{0.23\textwidth}
		\centering
		$\begin{bmatrix}
		* & - & + \\
		+ & 0 & + \\
		- & - & 0
		\end{bmatrix}$				
		\captionof{SP}{\textbf{DNA}}
	\end{minipage}
	\begin{minipage}{0.23\textwidth}
		\centering
		$\begin{bmatrix}
		- & - & - \\
		- & 0 & + \\
		+ & + & 0
		\end{bmatrix}$			
		\captionof{SP}{\textbf{DNA}}
	\end{minipage}
	
	Suppose that the digraph of $A=[a_{ij}]\in M_3(\mathbb{R})$ is Digraph \ref{D23}. If the sign pattern of $A$ is
	\begin{itemize}
		\item SP \ref{SP23.2}, then $A$ is AP if and only if $-a_{11}-\dfrac{a_{12}a_{23}}{a_{13}}, -a_{11}-\dfrac{a_{23}a_{31}}{a_{21}}, -\dfrac{a_{13}a_{21}}{a_{23}}, -a_{11}-\dfrac{a_{21}a_{32}}{a_{31}}, -\dfrac{a_{12}a_{31}}{a_{32}} < -a_{11}-\dfrac{a_{13}a_{32}}{a_{12}}$. In particular, if we take $\pm a_{11}=a_{13}=a_{21}=a_{23}=a_{31}=a_{32}=1$ and $a_{12}=-10$, then $A$ is not AP. On the other hand, if we take $\pm a_{11}=a_{13}=a_{21}=a_{23}=a_{31}=a_{32}=1$ and $a_{12}=-0.1$, then $A$ is AP.
		\item SP \ref{SP23.3}, then $A$ is AP if and only if $-a_{11}-\dfrac{a_{12}a_{23}}{a_{13}}, -a_{11}-\dfrac{a_{23}a_{31}}{a_{21}}, -a_{11}-\dfrac{a_{13}a_{32}}{a_{12}}, -a_{11}-\dfrac{a_{21}a_{32}}{a_{31}}, -\dfrac{a_{12}a_{31}}{a_{32}} < -\dfrac{a_{13}a_{21}}{a_{23}}$. In particular, if we take $\pm a_{11}=a_{12}=a_{13}=a_{21}=a_{31}=a_{32}=1$ and $a_{23}=-10$, then $A$ is not AP. On the other hand, if we take $\pm a_{11}=a_{12}=a_{13}=a_{21}=a_{31}=a_{32}=1$ and $a_{23}=-0.1$, then $A$ is AP.
		\item SP \ref{SP23.4}, then $A$ is AP if and only if $-a_{11}-\dfrac{a_{12}a_{23}}{a_{13}}, -\dfrac{a_{13}a_{21}}{a_{23}}, -a_{11}-\dfrac{a_{21}a_{32}}{a_{31}}, -\dfrac{a_{12}a_{31}}{a_{32}} < -a_{11}-\dfrac{a_{23}a_{31}}{a_{21}}, -a_{11}-\dfrac{a_{13}a_{32}}{a_{12}}$. In particular, if we take $\pm a_{11}=a_{13}=-a_{21}=a_{23}=a_{31}=a_{32}=1$ and $a_{12}=-10$, then $A$ is not AP. On the other hand, if we take $\pm a_{11}=a_{13}=a_{23}=a_{31}=a_{32}=1$ and $a_{12}=a_{21}=-0.1$, then $A$ is AP.
		\item SP \ref{SP23.5}, then $A$ is AP if and only if $-a_{11}-\dfrac{a_{12}a_{23}}{a_{13}}, -a_{11}-\dfrac{a_{23}a_{31}}{a_{21}}, -a_{11}-\dfrac{a_{21}a_{32}}{a_{31}}, -a_{11}-\dfrac{a_{13}a_{32}}{a_{12}} < -\dfrac{a_{13}a_{21}}{a_{23}}, -\dfrac{a_{12}a_{31}}{a_{32}}$. In particular, if we take $\pm a_{11}=a_{12}=a_{13}=a_{21}=a_{31}=-a_{32}=1$ and $a_{23}=-10$, then $A$ is not AP. On the other hand, if we take $\pm a_{11}=a_{12}=a_{13}=a_{21}=a_{31}=1$, $a_{23}=-0.1$, and $a_{32}=-0.3$, then $A$ is AP.
		\item SP \ref{SP23.6}, then $A$ is AP if and only if $-a_{11}-\dfrac{a_{12}a_{23}}{a_{13}}, -a_{11}-\dfrac{a_{23}a_{31}}{a_{21}}, -a_{11}-\dfrac{a_{21}a_{32}}{a_{31}}, -\dfrac{a_{12}a_{31}}{a_{32}} < -\dfrac{a_{13}a_{21}}{a_{23}}, -a_{11}-\dfrac{a_{13}a_{32}}{a_{12}}$. In particular, if we take $\pm a_{11}=a_{13}=a_{21}=a_{31}=a_{32}=1$ and $a_{12}=a_{23}=-10$, then $A$ is not AP. On the other hand, if we take $\pm a_{11}=a_{13}=a_{21}=a_{31}=a_{32}=1$, $a_{12}=-0.1$, and $a_{23}=-0.3$, then $A$ is AP.
		\item SP \ref{SP23.7}, then $A$ is AP if and only if $-a_{11}-\dfrac{a_{12}a_{23}}{a_{13}}, -a_{11}-\dfrac{a_{23}a_{31}}{a_{21}}, -\dfrac{a_{13}a_{21}}{a_{23}}, -\dfrac{a_{12}a_{31}}{a_{32}} < -a_{11}-\dfrac{a_{21}a_{32}}{a_{31}}, -a_{11}-\dfrac{a_{13}a_{32}}{a_{12}}$. In particular, if we take $\pm a_{11}=a_{13}=a_{21}=a_{23}=a_{32}=1$ and $a_{12}=a_{31}=-10$, then $A$ is not AP. On the other hand, if we take $\pm a_{11}=a_{13}=a_{21}=a_{23}=a_{32}=1$, $a_{12}=-0.1$, and $a_{31}=-0.3$, then $A$ is AP.
		\item SP \ref{SP23.8}, then $A$ is AP if and only if $-a_{11}-\dfrac{a_{12}a_{23}}{a_{13}}, -a_{11}-\dfrac{a_{23}a_{31}}{a_{21}},  -a_{11}-\dfrac{a_{13}a_{32}}{a_{12}} < -a_{11}-\dfrac{a_{21}a_{32}}{a_{31}}, -\dfrac{a_{13}a_{21}}{a_{23}}, -\dfrac{a_{12}a_{31}}{a_{32}}$. In particular, if we take $a_{11}=a_{12}=-a_{13}=a_{23}=a_{31}=a_{32}=1$ and $a_{21}=-10$, then $A$ is not AP. On the other hand, if we take $-a_{12}=-a_{13}=a_{23}=a_{31}=a_{32}=1$, $a_{11}=20$, and $a_{21}=-10$, then $A$ is AP.
		\item SP \ref{SP23.9}, then $A$ is AP if and only if $-a_{11}-\dfrac{a_{12}a_{23}}{a_{13}}, -a_{11}-\dfrac{a_{23}a_{31}}{a_{21}},  -\dfrac{a_{12}a_{31}}{a_{32}} < -a_{11}-\dfrac{a_{21}a_{32}}{a_{31}}, -\dfrac{a_{13}a_{21}}{a_{23}}, -a_{11}-\dfrac{a_{13}a_{32}}{a_{12}}$. In particular, if we take $-a_{12}=a_{13}=a_{21}=-a_{23}=-a_{31}=a_{32}1$ and $\pm a_{11}=10$, then $A$ is not AP. On the other hand, if we take $\pm a_{11}=-a_{12}=a_{13}=a_{21}=-a_{31}=a_{32}1$ and $a_{23}=10$, then $A$ is AP.
	\end{itemize}
	
	\item \textbf{Graphs with 8 directed edges}. There are two nonequivalent strongly connected digraphs with 8 directed edges. We list down the nonequivalent sign pattern matrices corresponding to each digraph and except for the AAP sign patterns, we use Theorems \ref{thm2}.2,\ref{thm2}.3 and \ref{thm4} to classify these sign patterns as indicated below. 
	
	\begin{minipage}{0.3\textwidth}
		\centering
		\includegraphics[scale=0.6]{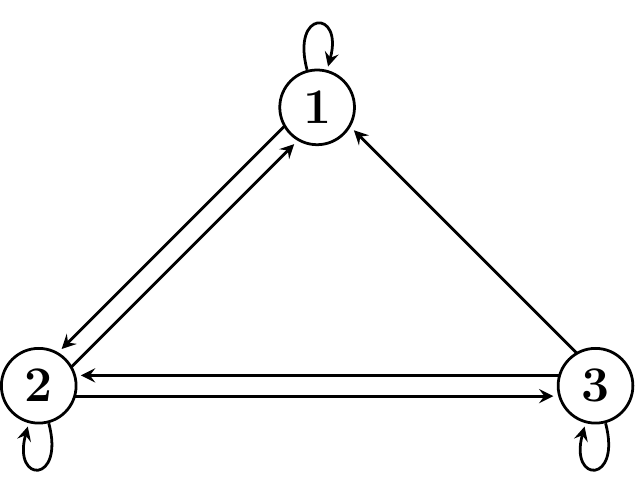}
		\captionof{Digraph}{}
		\label{D24}
	\end{minipage}
	\begin{minipage}{0.23\textwidth}
		\centering
		$\begin{bmatrix}
		* & + & 0 \\
		+ & * & + \\
		+ & + & *
		\end{bmatrix}$
		\captionof{SP}{\textbf{RAP}}
	\end{minipage}
	\begin{minipage}{0.23\textwidth}
		\centering
		$\begin{bmatrix}
		* & + & 0 \\
		- & * & + \\
		+ & + & *
		\end{bmatrix}$
		\captionof{SP}{\textbf{AAP}}
		\label{SP24.2}
	\end{minipage}
	\begin{minipage}{0.23\textwidth}
		\centering
		$\begin{bmatrix}
		* & + & 0 \\
		+ & * & + \\
		- & + & *
		\end{bmatrix}$
		\captionof{SP}{\textbf{AAP}}
		\label{SP24.3}
	\end{minipage}
	\smallskip\\
	\begin{minipage}{0.3\textwidth}
		\centering
		\ 
	\end{minipage}
	\begin{minipage}{0.23\textwidth}
		\centering
		$\begin{bmatrix}
		- & - & 0 \\
		- & * & + \\
		+ & + & -
		\end{bmatrix}$
		\captionof{SP}{\textbf{AAP}}
		\label{SP24.4}
	\end{minipage}
	\begin{minipage}{0.23\textwidth}
		\centering
		$\begin{bmatrix}
		+ & - & 0 \\
		- & * & + \\
		+ & + & *
		\end{bmatrix}$				
		\captionof{SP}{\textbf{AAP}}
		\label{SP24.5}
	\end{minipage}
	\begin{minipage}{0.23\textwidth}
		\centering
		$\begin{bmatrix}
		- & - & 0 \\
		+ & - & + \\
		- & + & -
		\end{bmatrix}$			
		\captionof{SP}{\textbf{AAP}}
		\label{SP24.6}
	\end{minipage}
	\smallskip\\
	\begin{minipage}{0.3\textwidth}
		\centering
		\ 
	\end{minipage}
	\begin{minipage}{0.23\textwidth}
		\centering
		$\begin{bmatrix}
		+ & - & 0 \\
		+ & * & + \\
		- & + & *
		\end{bmatrix}$
		\captionof{SP}{\textbf{AAP}}
		\label{SP24.7}
	\end{minipage}
	\begin{minipage}{0.23\textwidth}
		\centering
		$\begin{bmatrix}
		* & + & 0 \\
		- & * & + \\
		+ & - & *
		\end{bmatrix}$				
		\captionof{SP}{\textbf{AAP}}
		\label{SP24.8}
	\end{minipage}
	\begin{minipage}{0.23\textwidth}
		\centering
		$\begin{bmatrix}
		* & - & 0 \\
		+ & * & + \\
		+ & + & *
		\end{bmatrix}$			
		\captionof{SP}{\textbf{DNA}}
	\end{minipage}
	\begin{minipage}{0.3\textwidth}
		\centering
		\ 
	\end{minipage}
	\begin{minipage}{0.23\textwidth}
		\centering
		$\begin{bmatrix}
		* & - & 0 \\
		+ & * & - \\
		+ & + & *
		\end{bmatrix}$
		\captionof{SP}{\textbf{DNA}}
	\end{minipage}
	\begin{minipage}{0.23\textwidth}
		\centering
		$\begin{bmatrix}
		* & - & 0 \\
		+ & * & + \\
		+ & - & *
		\end{bmatrix}$				
		\captionof{SP}{\textbf{DNA}}
	\end{minipage}
	\begin{minipage}{0.23\textwidth}
		\centering
		$\begin{bmatrix}
		* & + & 0 \\
		- & * & + \\
		- & + & *
		\end{bmatrix}$			
		\captionof{SP}{\textbf{DNA}}
	\end{minipage}
	\begin{minipage}{0.3\textwidth}
		\centering
		\ 
	\end{minipage}
	\begin{minipage}{0.23\textwidth}
		\centering
		$\begin{bmatrix}
		- & - & 0 \\
		- & * & + \\
		+ & + & +
		\end{bmatrix}$
		\captionof{SP}{\textbf{DNA}}
	\end{minipage}
	\begin{minipage}{0.23\textwidth}
		\centering
		$\begin{bmatrix}
		- & - & 0 \\
		+ & - & + \\
		- & + & +
		\end{bmatrix}$				
		\captionof{SP}{\textbf{DNA}}
	\end{minipage}
	\begin{minipage}{0.23\textwidth}
		\centering
		$\begin{bmatrix}
		- & - & 0 \\
		+ & + & + \\
		- & + & *
		\end{bmatrix}$			
		\captionof{SP}{\textbf{DNA}}
	\end{minipage}
	
	Suppose that the digraph of $A=[a_{ij}]\in M_3(\mathbb{R})$ is Digraph \ref{D24}. If the sign pattern of $A$ is
	\begin{itemize}
		\item SP \ref{SP24.2}, then $A$ is AP if and only if $-a_{11}-a{22}, -a{22}-a_{33}, -a_{11}-a_{33}-\dfrac{a_{21}a_{32}}{a_{31}}, -a{22}-a_{33}-\dfrac{a_{12}a_{31}}{a_{32}} < -a_{11}-a{22}-\dfrac{a_{23}a_{31}}{a_{21}}$. In particular, if we take $\pm a_{11}=a_{12}=\pm a_{22}=a_{23}=a_{31}=1$ and $-a_{21}=a_{32}=\pm a_{33}=10$, then $A$ is not AP. On the other hand, if we take $\pm a_{11}=a_{12}=\pm a_{22}=a_{23}=a_{31}=a_{32}=\pm a_{33}=1$ and $a_{21}=-0.1$, then $A$ is AP.
		\item SP \ref{SP24.3}, then $A$ is AP if and only if $-a_{11}-a{22}, -a{22}-a_{33}, -a_{11}-a{22}-\dfrac{a_{23}a_{31}}{a_{21}}, -a{22}-a_{33}-\dfrac{a_{12}a_{31}}{a_{32}} < -a_{11}-a_{33}-\dfrac{a_{21}a_{32}}{a_{31}}$. In particular, if we take $\pm a_{11}=a_{21}=\pm a_{22}=a_{23}=a_{32}=\pm a_{33}=1$ and $a_{12}=-a_{31}=10$, then $A$ is not AP. On the other hand, if we take $\pm a_{11}=a_{12}=a_{21}=\pm a_{22}=a_{23}=a_{32}=\pm a_{33}=1$ and $a_{31}=-0.1$, then $A$ is AP.
		\item SP \ref{SP24.4}, then $A$ is AP if and only if $-a_{11}-a{22}, -a_{11}-a{22}-\dfrac{a_{23}a_{31}}{a_{21}} < -a{22}-a_{33}, -a_{11}-a_{33}-\dfrac{a_{21}a_{32}}{a_{31}}, -a{22}-a_{33}-\dfrac{a_{12}a_{31}}{a_{32}}$. In particular, if we take $-a_{11}=-a_{12}=-a_{21}=\pm a_{22}=a_{23}=a_{31}=a_{32}=-a_{33}=1$, then $A$ is not AP. On the other hand, if we take $-a_{11}=-a_{12}=\pm a_{22}=a_{23}=a_{31}=a_{32}=1$, $a_{21}=-10$, and $a_{33}=-2$, then $A$ is AP.
		\item SP \ref{SP24.5}, then $A$ is AP if and only if $-a_{11}-a{22}, -a_{11}-a{22}-\dfrac{a_{23}a_{31}}{a_{21}} < -a{22}-a_{33}, -a_{11}-a_{33}-\dfrac{a_{21}a_{32}}{a_{31}}, -a{22}-a_{33}-\dfrac{a_{12}a_{31}}{a_{32}}$. In particular, if we take $a_{11}=-a_{12}=\pm a_{22}=a_{23}=a_{31}=a_{32}=\pm a_{33}=1$ and $a_{21}=-0.1$, then $A$ is not AP. On the other hand, if we take $-a_{12}=\pm a_{22}=a_{23}=a_{31}=a_{32}=\pm a_{33}=1$, $a_{21}=-10$, and $a_{11}=2$, then $A$ is AP.
		\item SP \ref{SP24.6}, then $A$ is AP if and only if $-a_{11}-a{22}, -a_{11}-a_{33}-\dfrac{a_{21}a_{32}}{a_{31}} < -a{22}-a_{33}, -a_{11}-a{22}-\dfrac{a_{23}a_{31}}{a_{21}}, -a{22}-a_{33}-\dfrac{a_{12}a_{31}}{a_{32}}$. In particular, if we take $-a_{11}=-a_{12}=a_{21}=-a_{22}=a_{23}=-a_{31}=a_{32}=-a_{33}=1$, then $A$ is not AP. On the other hand, if we take $-a_{11}=-a_{12}=a_{21}=a_{23}=1$, $a_{22}=a_{31}=-10$, $a_{32}=12$, and $a_{33}=2$, then $A$ is AP.
		\item SP \ref{SP24.7}, then $A$ is AP if and only if $-a_{11}-a{22}, -a_{11}-a_{33}-\dfrac{a_{21}a_{32}}{a_{31}} < -a{22}-a_{33}, -a_{11}-a{22}-\dfrac{a_{23}a_{31}}{a_{21}}, -a{22}-a_{33}-\dfrac{a_{12}a_{31}}{a_{32}}$. In particular, if we take $a_{11}=-a_{12}=a_{21}=\pm a_{22}=a_{23}=-a_{31}=a_{32}=\pm a_{33}=1$, then $A$ is not AP. On the other hand, if we take $a_{21}=\pm a_{22}=a_{23}=a_{32}=\pm a_{33}=1$, $a_{31}=-10$, $a_{12}=-0.1$, and $a_{11}=3$, then $A$ is AP.
		\item SP \ref{SP24.8}, then $A$ is AP if and only if $-a_{11}-a{22}, -a{22}-a_{33}, -a_{11}-a_{33}-\dfrac{a_{21}a_{32}}{a_{31}} < -a_{11}-a{22}-\dfrac{a_{23}a_{31}}{a_{21}}, -a{22}-a_{33}-\dfrac{a_{12}a_{31}}{a_{32}}$. In particular, if we take $\pm a_{11}=a_{12}=-a_{21}=a_{23}=a_{31}=1$, $\pm a_{22}=-a_{32}=10$, and $\pm a_{33}=20$, then $A$ is not AP. On the other hand, if we take $\pm a_{11}=a_{12}=\pm a_{22}, a_{23}=a_{31}=\pm a_{33}=1$ and $a_{21}=a_{32}=-0.1$, then $A$ is AP.
	\end{itemize}
	
	\begin{minipage}{0.3\textwidth}
		\centering
		\includegraphics[scale=0.6]{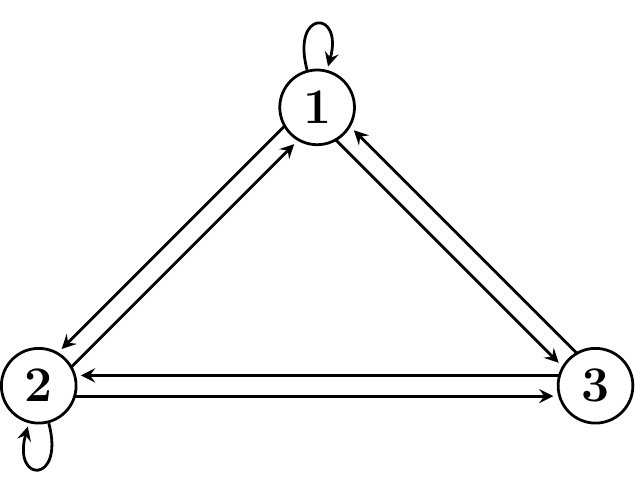}
		\captionof{Digraph}{}
		\label{D25}
	\end{minipage}
	\begin{minipage}{0.23\textwidth}
		\centering
		$\begin{bmatrix}
		* & + & + \\
		+ & * & + \\
		+ & + & 0
		\end{bmatrix}$
		\captionof{SP}{\textbf{RAP}}
	\end{minipage}
	\begin{minipage}{0.23\textwidth}
		\centering
		$\begin{bmatrix}
		* & - & + \\
		+ & * & + \\
		+ & + & 0
		\end{bmatrix}$
		\captionof{SP}{\textbf{AAP}}
		\label{SP25.2}
	\end{minipage}
	\begin{minipage}{0.23\textwidth}
		\centering
		$\begin{bmatrix}
		* & + & - \\
		+ & * & + \\
		+ & + & 0
		\end{bmatrix}$
		\captionof{SP}{\textbf{AAP}}
		\label{SP25.3}
	\end{minipage}
	\smallskip\\
	\begin{minipage}{0.3\textwidth}
		\centering
		\ 
	\end{minipage}
	\begin{minipage}{0.23\textwidth}
		\centering
		$\begin{bmatrix}
		* & - & + \\
		- & * & + \\
		+ & + & 0
		\end{bmatrix}$
		\captionof{SP}{\textbf{AAP}}
		\label{SP25.4}
	\end{minipage}
	\begin{minipage}{0.23\textwidth}
		\centering
		$\begin{bmatrix}
		* & + & - \\
		+ & * & + \\
		- & + & 0
		\end{bmatrix}$				
		\captionof{SP}{\textbf{AAP}}
		\label{SP25.5}
	\end{minipage}
	\begin{minipage}{0.23\textwidth}
		\centering
		$\begin{bmatrix}
		* & - & + \\
		+ & * & - \\
		+ & + & 0
		\end{bmatrix}$			
		\captionof{SP}{\textbf{AAP}}
		\label{SP25.6}
	\end{minipage}
	\smallskip\\
	\begin{minipage}{0.3\textwidth}
		\centering
		\ 
	\end{minipage}
	\begin{minipage}{0.23\textwidth}
		\centering
		$\begin{bmatrix}
		* & + & + \\
		+ & * & - \\
		- & + & 0
		\end{bmatrix}$
		\captionof{SP}{\textbf{AAP}}
		\label{SP25.7}
	\end{minipage}
	\begin{minipage}{0.23\textwidth}
		\centering
		$\begin{bmatrix}
		* & - & + \\
		+ & * & - \\
		- & + & 0
		\end{bmatrix}$				
		\captionof{SP}{\textbf{AAP}}
		\label{SP25.8}
	\end{minipage}
	\begin{minipage}{0.23\textwidth}
		\centering
		$\begin{bmatrix}
		+ & - & - \\
		- & * & + \\
		+ & + & 0
		\end{bmatrix}$			
		\captionof{SP}{\textbf{AAP}}
		\label{SP25.9}
	\end{minipage}
	\begin{minipage}{0.3\textwidth}
		\centering
		\ 
	\end{minipage}
	\begin{minipage}{0.23\textwidth}
		\centering
		$\begin{bmatrix}
		+ & - & - \\
		+ & * & + \\
		- & + & 0
		\end{bmatrix}$
		\captionof{SP}{\textbf{AAP}}
		\label{SP25.10}
	\end{minipage}
	\begin{minipage}{0.23\textwidth}
		\centering
		$\begin{bmatrix}
		* & - & - \\
		+ & * & + \\
		+ & + & 0
		\end{bmatrix}$				
		\captionof{SP}{\textbf{DNA}}
	\end{minipage}
	\begin{minipage}{0.23\textwidth}
		\centering
		$\begin{bmatrix}
		* & + & + \\
		+ & * & + \\
		- & - & 0
		\end{bmatrix}$			
		\captionof{SP}{\textbf{DNA}}
	\end{minipage}
	\begin{minipage}{0.3\textwidth}
		\centering
		\ 
	\end{minipage}
	\begin{minipage}{0.23\textwidth}
		\centering
		$\begin{bmatrix}
		* & - & - \\
		+ & * & - \\
		+ & + & 0
		\end{bmatrix}$
		\captionof{SP}{\textbf{DNA}}
	\end{minipage}
	\begin{minipage}{0.23\textwidth}
		\centering
		$\begin{bmatrix}
		* & - & - \\
		+ & * & + \\
		+ & - & 0
		\end{bmatrix}$				
		\captionof{SP}{\textbf{DNA}}
	\end{minipage}
	\begin{minipage}{0.23\textwidth}
		\centering
		$\begin{bmatrix}
		- & - & - \\
		- & * & + \\
		+ & + & 0
		\end{bmatrix}$			
		\captionof{SP}{\textbf{DNA}}
	\end{minipage}
	\begin{minipage}{0.76\textwidth}
		\centering
		\ 
	\end{minipage}
	\begin{minipage}{0.23\textwidth}
		\centering
		$\begin{bmatrix}
		- & - & - \\
		+ & * & + \\
		- & + & 0
		\end{bmatrix}$			
		\captionof{SP}{\textbf{DNA}}
	\end{minipage}
	
	Suppose that the digraph of $A=[a_{ij}]\in M_3(\mathbb{R})$ is Digraph \ref{D25}. If the sign pattern of $A$ is
	\begin{itemize}
		\item SP \ref{SP25.2}, then $A$ is AP if and only if $-a_{11}-a_{22}-\dfrac{a_{23}a_{31}}{a_{21}}, -a_{11}-\dfrac{a_{21}a_{32}}{a_{31}}, -a_{22}-\dfrac{a_{12}a_{31}}{a_{32}}, -a_{11}-\dfrac{a_{12}a_{23}}{a_{13}}, -a_{22}-\dfrac{a_{13}a_{21}}{a_{23}} < -a_{11}-a_{22}-\dfrac{a_{13}a_{32}}{a_{12}}$. In particular, if we take $\pm a_{11}=a_{13}=a_{21}=\pm a_{22}=a_{23}=a_{31}=a_{32}=1$ and $a_{12}=-10$, then $A$ is not AP. On the other hand, if we take $\pm a_{11}=a_{13}=a_{21}=\pm a_{22}=a_{23}=a_{31}=a_{32}=1$ and $a_{12}=-0.1$, then $A$ is AP.
		\item SP \ref{SP25.3}, then $A$ is AP if and only if $-a_{11}-a_{22}-\dfrac{a_{23}a_{31}}{a_{21}}, -a_{11}-\dfrac{a_{21}a_{32}}{a_{31}}, -a_{22}-\dfrac{a_{12}a_{31}}{a_{32}}, -a_{11}-a_{22}-\dfrac{a_{13}a_{32}}{a_{12}}, -a_{22}-\dfrac{a_{13}a_{21}}{a_{23}} < -a_{11}-\dfrac{a_{12}a_{23}}{a_{13}}$. In particular, if we take $\pm a_{11}=a_{12}=a_{21}=\pm a_{22}=a_{23}=a_{31}=a_{32}=1$ and $a_{13}=-10$, then $A$ is not AP. On the other hand, if we take $\pm a_{11}=a_{13}=a_{21}=\pm a_{22}=a_{23}=a_{31}=a_{32}=1$ and $a_{13}=-0.1$, then $A$ is AP.
		\item SP \ref{SP25.4}, then $A$ is AP if and only if $-a_{11}-\dfrac{a_{12}a_{23}}{a_{13}}, -a_{11}-\dfrac{a_{21}a_{32}}{a_{31}}, -a_{22}-\dfrac{a_{12}a_{31}}{a_{32}},  -a_{22}-\dfrac{a_{13}a_{21}}{a_{23}} < -a_{11}-a_{22}-\dfrac{a_{13}a_{32}}{a_{12}}, -a_{11}-a_{22}-\dfrac{a_{23}a_{31}}{a_{21}}$. In particular, if we take $\pm a_{11}=-a_{12}=a_{13}=-a_{21}=\pm a_{22}=a_{31}=a_{32}=1$ and $a_{23}=2$, then $A$ is not AP. On the other hand, if we take $\pm a_{11}=a_{13}=\pm a_{22}=a_{23}=a_{31}=a_{32}=1$ and $a_{12}=a_{21}=-0.1$, then $A$ is AP.
		\item SP \ref{SP25.5}, then $A$ is AP if and only if $-a_{11}-a_{22}-\dfrac{a_{13}a_{32}}{a_{12}}, -a_{11}-a_{22}-\dfrac{a_{23}a_{31}}{a_{21}}, -a_{22}-\dfrac{a_{12}a_{31}}{a_{32}},  -a_{22}-\dfrac{a_{13}a_{21}}{a_{23}} < -a_{11}-\dfrac{a_{12}a_{23}}{a_{13}}, -a_{11}-\dfrac{a_{21}a_{32}}{a_{31}}$. In particular, if we take $\pm a_{11}=a_{12}=-a_{13}=a_{21}=\pm a_{22}=-a_{31}=a_{32}=1$ and $a_{23}=2$, then $A$ is not AP. On the other hand, if we take $\pm a_{11}=a_{12}=a_{21}=\pm a_{22}=a_{23}=a_{32}=1$ and $a_{13}=a_{31}=-0.1$, then $A$ is AP.
		\item SP \ref{SP25.6}, then $A$ is AP if and only if $-a_{11}-\dfrac{a_{12}a_{23}}{a_{13}}, -a_{11}-\dfrac{a_{21}a_{32}}{a_{31}}, -a_{11}-a_{22}-\dfrac{a_{23}a_{31}}{a_{21}}, -a_{22}-\dfrac{a_{12}a_{31}}{a_{32}} < -a_{11}-a_{22}-\dfrac{a_{13}a_{32}}{a_{12}}, -a_{22}-\dfrac{a_{13}a_{21}}{a_{23}}$. In particular, if we take $\pm a_{11}=-a_{12}=a_{13}=a_{21}=\pm a_{22}=-a_{23}=a_{32}=1$ and $a_{31}=10$, then $A$ is not AP. On the other hand, if we take $\pm a_{11}=a_{13}=a_{21}=\pm a_{22}=a_{31}=a_{32}=1$ and $a_{12}=a_{23}=-0.1$, then $A$ is AP.
		\item SP \ref{SP25.7}, then $A$ is AP if and only if $-a_{11}-\dfrac{a_{12}a_{23}}{a_{13}}, -a_{11}-a_{22}-\dfrac{a_{13}a_{32}}{a_{12}}, -a_{11}-a_{22}-\dfrac{a_{23}a_{31}}{a_{21}}, -a_{22}-\dfrac{a_{12}a_{31}}{a_{32}} < -a_{11}-\dfrac{a_{21}a_{32}}{a_{31}}, -a_{22}-\dfrac{a_{13}a_{21}}{a_{23}}$. In particular, if we take $\pm a_{11}=a_{13}=a_{21}=\pm a_{22}=-a_{23}=-a_{31}=a_{32}=1$ and $a_{12}=2$, then $A$ is not AP. On the other hand, if we take $\pm a_{11}=a_{12}=a_{13}=a_{21}=\pm a_{22}=a_{32}=1$ and $a_{23}=a_{31}=-0.1$, then $A$ is AP.
		\item SP \ref{SP25.8}, then $A$ is AP if and only if $-a_{11}-\dfrac{a_{12}a_{23}}{a_{13}}, -a_{11}-a_{22}-\dfrac{a_{23}a_{31}}{a_{21}}, -a_{22}-\dfrac{a_{12}a_{31}}{a_{32}} < -a_{11}-a_{22}-\dfrac{a_{13}a_{32}}{a_{12}}, -a_{11}-\dfrac{a_{21}a_{32}}{a_{31}}, -a_{22}-\dfrac{a_{13}a_{21}}{a_{23}}$. In particular, if we take $-a_{12}=a_{13}=a_{21}=-a_{23}=-a_{31}=a_{32}=1$, $\pm a_{11}=5$, and $\pm a_{22}=10$, then $A$ is not AP. On the other hand, if we take $\pm a_{11}=a_{13}=a_{21}=\pm a_{22}=a_{32}=1$ and $a_{12}=a_{23}=a_{31}=-0.1$, then $A$ is AP.
		\item SP \ref{SP25.9}, then $A$ is AP if and only if $-a_{11}-\dfrac{a_{21}a_{32}}{a_{31}}, -a_{22}-\dfrac{a_{13}a_{21}}{a_{23}}, -a_{22}-\dfrac{a_{12}a_{31}}{a_{32}} < -a_{11}-a_{22}-\dfrac{a_{13}a_{32}}{a_{12}}, -a_{11}-\dfrac{a_{12}a_{23}}{a_{13}}, -a_{11}-a_{22}-\dfrac{a_{23}a_{31}}{a_{21}}$. In particular, if we take $a_{11}=-a_{12}=-a_{13}=-a_{21}=\pm a_{22}=a_{23}=a_{31}=a_{32}=1$, then $A$ is not AP. On the other hand, if we take $-a_{12}=-a_{13}=\pm a_{22}=a_{23}=a_{31}=a_{32}=1$, $a_{11}=10$, and $a_{21}=-5$, then $A$ is AP.
		\item SP \ref{SP25.10}, then $A$ is AP if and only if $-a_{11}-a_{22}-\dfrac{a_{23}a_{31}}{a_{21}}, -a_{22}-\dfrac{a_{13}a_{21}}{a_{23}}, -a_{22}-\dfrac{a_{12}a_{31}}{a_{32}} < -a_{11}-a_{22}-\dfrac{a_{13}a_{32}}{a_{12}}, -a_{11}-\dfrac{a_{12}a_{23}}{a_{13}}, -a_{11}-\dfrac{a_{21}a_{32}}{a_{31}}$. In particular, if we take $a_{11}=-a_{12}=-a_{13}=a_{21}=\pm a_{22}=a_{23}=-a_{31}=a_{32}=1$, then $A$ is not AP. On the other hand, if we take $-a_{12}=-a_{13}=a_{21}=\pm a_{22}=a_{23}=a_{32}=1$, $a_{11}=10$, and $a_{31}=-2$, then $A$ is AP.
	\end{itemize}
	
	\item \textbf{Graphs with 9 directed edges}. There is only one nonequivalent strongly connected digraph with 9 directed edges. We list down the nonequivalent sign pattern matrices corresponding to the digraph and except for the AAP sign patterns, we use Theorems \ref{thm2}.2, \ref{thm2}.3, and \ref{thm4} to classify these sign patterns as indicated below.
	
	\begin{minipage}{0.3\textwidth}
		\centering
		\includegraphics[scale=0.6]{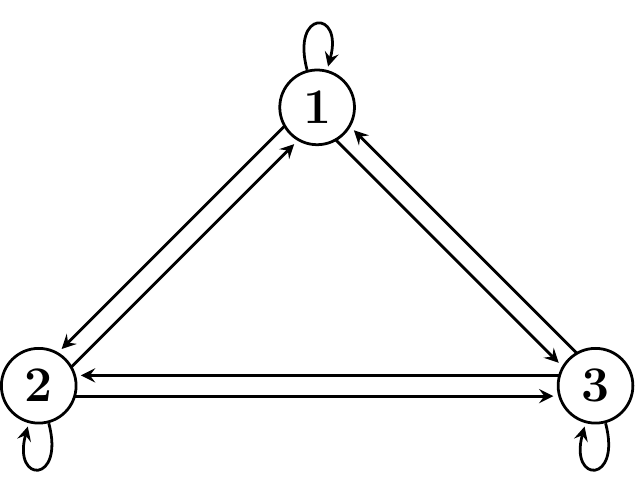}
		\captionof{Digraph}{}
		\label{D26}
	\end{minipage}
	\begin{minipage}{0.23\textwidth}
		\centering
		$\begin{bmatrix}
		* & + & + \\
		+ & * & + \\
		+ & + & *
		\end{bmatrix}$
		\captionof{SP}{\textbf{RAP}}
	\end{minipage}
	\begin{minipage}{0.23\textwidth}
		\centering
		$\begin{bmatrix}
		* & - & + \\
		+ & * & + \\
		+ & + & *
		\end{bmatrix}$
		\captionof{SP}{\textbf{AAP}}
		\label{SP26.2}
	\end{minipage}
	\begin{minipage}{0.23\textwidth}
		\centering
		$\begin{bmatrix}
		* & - & + \\
		+ & * & + \\
		+ & - & *
		\end{bmatrix}$
		\captionof{SP}{\textbf{DNA}}
	\end{minipage}
	\smallskip\\
	\begin{minipage}{0.3\textwidth}
		\centering
		\ 
	\end{minipage}
	\begin{minipage}{0.23\textwidth}
		\centering
		$\begin{bmatrix}
		* & - & + \\
		- & * & + \\
		+ & + & *
		\end{bmatrix}$
		\captionof{SP}{\textbf{AAP}}
		\label{SP26.4}
	\end{minipage}
	\begin{minipage}{0.23\textwidth}
		\centering
		$\begin{bmatrix}
		* & - & + \\
		+ & * & - \\
		+ & + & *
		\end{bmatrix}$				
		\captionof{SP}{\textbf{AAP}}
		\label{SP26.5}
	\end{minipage}
	\begin{minipage}{0.23\textwidth}
		\centering
		$\begin{bmatrix}
		* & - & - \\
		+ & * & - \\
		+ & + & *
		\end{bmatrix}$			
		\captionof{SP}{\textbf{DNA}}
	\end{minipage}
	\smallskip\\
	\begin{minipage}{0.3\textwidth}
		\centering
		\ 
	\end{minipage}
	\begin{minipage}{0.23\textwidth}
		\centering
		$\begin{bmatrix}
		* & - & + \\
		+ & * & - \\
		- & + & *
		\end{bmatrix}$
		\captionof{SP}{\textbf{AAP}}
		\label{SP26.7}
	\end{minipage}
	\begin{minipage}{0.23\textwidth}
		\centering
		$\begin{bmatrix}
		* & - & + \\
		- & * & - \\
		+ & + & -
		\end{bmatrix}$				
		\captionof{SP}{\textbf{AAP}}
		\label{SP26.8}
	\end{minipage}
	\begin{minipage}{0.23\textwidth}
		\centering
		$\begin{bmatrix}
		* & - & + \\
		- & - & - \\
		+ & + & +
		\end{bmatrix}$			
		\captionof{SP}{\textbf{DNA}}
	\end{minipage}
	
	Suppose that the digraph of $A=[a_{ij}]\in M_3(\mathbb{R})$ is Digraph \ref{D26}. If the sign pattern of $A$ is
	\begin{itemize}
		\item SP \ref{SP26.2}, then $A$ is AP if and only if $-a_{11}-a_{33}-\dfrac{a_{12}a_{23}}{a_{13}}, -a_{11}-a_{22}-\dfrac{a_{23}a_{31}}{a_{21}}, -a_{22}-a_{33}-\dfrac{a_{13}a_{21}}{a_{23}}, -a_{11}-a_{33}-\dfrac{a_{21}a_{32}}{a_{31}}, -a_{22}-a_{33}-\dfrac{a_{12}a_{31}}{a_{32}} < -a_{11}-a_{22}-\dfrac{a_{13}a_{32}}{a_{12}}$. In particular, if we take $\pm a_{11}=a_{13}=a_{21}=\pm a_{22}=a_{23}=a_{31}=a_{32}=\pm a_{33}=1$ and $a_{12}=-10$, then $A$ is not AP. On the other hand, if we take $\pm a_{11}=a_{13}=a_{21}=\pm a_{22}=a_{23}=a_{31}=a_{32}=\pm a_{33}=1$ and $a_{12}=-0.1$, then $A$ is AP.
		\item SP \ref{SP26.4}, then $A$ is AP if and only if $-a_{11}-a_{33}-\dfrac{a_{12}a_{23}}{a_{13}}, -a_{22}-a_{33}-\dfrac{a_{13}a_{21}}{a_{23}}, -a_{11}-a_{33}-\dfrac{a_{21}a_{32}}{a_{31}}, -a_{22}-a_{33}-\dfrac{a_{12}a_{31}}{a_{32}} < -a_{11}-a_{22}-\dfrac{a_{13}a_{32}}{a_{12}}, -a_{11}-a_{22}-\dfrac{a_{23}a_{31}}{a_{21}}$. In particular, if we take $\pm a_{11}=-a_{12}=\pm a_{22}=a_{23}=a_{31}=a_{32}=\pm a_{33}=1$ and $a_{13}=-a_{21}=10$, then $A$ is not AP. On the other hand, if we take $\pm a_{11}=a_{13}=\pm a_{22}=a_{23}=a_{31}=a_{32}=\pm a_{33}=1$ and $a_{12}=a_{21}=-0.1$, then $A$ is AP.
		\item SP \ref{SP26.5}, then $A$ is AP if and only if $-a_{11}-a_{33}-\dfrac{a_{12}a_{23}}{a_{13}}, -a_{11}-a_{22}-\dfrac{a_{23}a_{31}}{a_{21}}, -a_{11}-a_{33}-\dfrac{a_{21}a_{32}}{a_{31}}, -a_{22}-a_{33}-\dfrac{a_{12}a_{31}}{a_{32}} < -a_{11}-a_{22}-\dfrac{a_{13}a_{32}}{a_{12}}, -a_{22}-a_{33}-\dfrac{a_{13}a_{21}}{a_{23}}$. In particular, if we take $\pm a_{11}=-a_{12}=a_{13}=a_{21}=\pm a_{22}=-a_{23}=a_{31}=a_{32}=\pm a_{33}=1$, then $A$ is not AP. On the other hand, if we take $\pm a_{11}=a_{13}=a_{21}=\pm a_{22}=a_{31}=a_{32}=\pm a_{33}=1$ and $a_{12}=a_{23}=-0.1$, then $A$ is AP.
		\item SP \ref{SP26.7}, then $A$ is AP if and only if $-a_{11}-a_{33}-\dfrac{a_{12}a_{23}}{a_{13}}, -a_{11}-a_{22}-\dfrac{a_{23}a_{31}}{a_{21}}, -a_{22}-a_{33}-\dfrac{a_{12}a_{31}}{a_{32}} < -a_{11}-a_{22}-\dfrac{a_{13}a_{32}}{a_{12}}, -a_{22}-a_{33}-\dfrac{a_{13}a_{21}}{a_{23}}, -a_{11}-a_{33}-\dfrac{a_{21}a_{32}}{a_{31}}$. In particular, if we take $\pm a_{11}=-a_{12}=a_{13}=a_{21}=-a_{23}=-a_{31}=a_{32}=1$, $\pm a_{11}=1$, and $\pm a_{33}=10$, then $A$ is not AP. On the other hand, if we take $\pm a_{11}=a_{13}=a_{21}=\pm a_{22}=a_{32}=\pm a_{33}=1$ and $a_{12}=a_{23}=a_{31}=-0.1$, then $A$ is AP.
		\item SP \ref{SP26.8}, then $A$ is AP if and only if $-a_{11}-a_{22}-\dfrac{a_{13}a_{32}}{a_{12}}, -a_{22}-a_{33}-\dfrac{a_{13}a_{21}}{a_{23}}, -a_{11}-a_{22}-\dfrac{a_{23}a_{31}}{a_{21}} < -a_{11}-a_{33}-\dfrac{a_{12}a_{23}}{a_{13}},  -a_{11}-a_{33}-\dfrac{a_{21}a_{32}}{a_{31}}, -a_{22}-a_{33}-\dfrac{a_{12}a_{31}}{a_{32}}$. In particular, if we take $\pm a_{11}=-a_{12}=a_{13}=-a_{21}=\pm a_{22}=-a_{23}=a_{31}=a_{32}=-a_{33}=1$,then $A$ is not AP. On the other hand, if we take $\pm a_{11}=-a_{12}=-a_{21}=\pm a_{22}=-a_{23}=a_{31}=a_{32}=1$, $a_{13}=5$, and $a_{33}=-10$, then $A$ is AP.
\end{itemize}
	
\end{enumerate}

\begin{thebibliography}{WWW}

\bibitem{Kirk1}
S.~{Kirkland}, P.~{Qiao}, and X.~{Zhan}, Algebraically positive matrices, {\em Linear Algebra Appl.}, vol.~504, pp.~14--26, 2016.

\bibitem{Pillai}
S.U.~{Pillai}, T.~{Suel}, and S.~{Cha}, The Perron-Frobenius theorem: some of its applications, {\em IEEE Signal Processing Magazine}, vol.~22, no.~2, pp.~62--75, 2005.

\bibitem{MacCluer}
C.R.~{MacCluer}, The Many Proofs and Applications of Perron's Theorem, {\em SIAM Review}, vol.~42, no.~3, pp.~487--498, 2000.


\bibitem{handbook}
L.~{Hogben}, ``Handbook of Linear Algebra Second Edition,'' {\em CRC Press}, 2014.

\end{thebibliography}
\end{document}